\documentclass
{amsart}
\usepackage[cmtip,all]{xy}
\usepackage{graphicx}
\usepackage{amssymb}

\usepackage{graphics}
\usepackage{epsf,epsfig,amsmath}
\usepackage{psfrag}

\newtheorem{theorem}{Theorem}[section]
\newtheorem{lemma}[theorem]{Lemma}
\newtheorem{prop}[theorem]{Proposition}

\newtheorem{cor}[theorem]{Corollary}
\theoremstyle{definition}
\newtheorem{definition}[theorem]{Definition}
\newtheorem{example}[theorem]{Example}

\theoremstyle{remark}
\newtheorem{remark}[theorem]{Remark}

\numberwithin{equation}{section}

\newcommand{\cF}{\mathcal F}

\newcommand{\ics}{x_1 , \ldots , x_n}
\newcommand{\icss}{x_1^2 , \ldots , x_n^2}
\newcommand{\lcm}{{\rm lcm}}

\newcommand{\Rn}{\mathbb{Z}[X]}

\begin{document}

\title{Cores of simplicial complexes}

\author{Mario Marietti}
\address{Universit\`a degli Studi di Roma ``La Sapienza'', Piazzale A. Moro 5, 00185 Roma}
\email{marietti@mat.uniroma1.it \hspace{10pt} www.mat.uniroma1.it/$\sim$marietti}

\author{Damiano Testa}
\address{Universit\`a degli Studi di Roma ``La Sapienza'', Piazzale A. Moro 5, 00185 Roma}
\email{testa@mat.uniroma1.it \hspace{23pt} www.mat.uniroma1.it/$\sim$testa}

\begin{abstract}
We introduce a method to reduce the study of the topology of a simplicial complex to that of a simpler one. We give some applications of this method to complexes arising from graphs. 
As a consequence, we answer some questions raised in  [European J.~Combin.~27 (2006), no.~6, 906-923] on the independence complex and the dominance complex.  The techniques used come mainly from monomial ideal theory.
\end{abstract}


\maketitle

\section{Introduction}

In this paper we study the topology of a simplicial complex $\Delta$ by introducing a family 
$c(\Delta)$ of complexes that we call the {\it core of $\Delta $}.  The homotopy types of 
$\Delta $ and its cores are closely related.  Indeed we show that one of the following 
happens: either $\Delta $ collapses onto a point, or it is simple-homotopic to an 
iterated suspension of any element of its core.  Since we work with the Stanley-Reisner 
ideal associated to $\Delta$, the techniques used come mainly from monomial ideal theory.

The application motivating this method is in the study of the independence and dominance 
complexes of a graph $G$ and allows us to answer some questions posed in~\cite{EH}, where 
Ehrenborg and Hetyei prove what follows:
\begin{enumerate}
\item[-] the independence complex of a forest $F$ is always contractible or homotopic to a sphere;        
\item[-] the dominance complex of a forest $F$ is always homotopic to a sphere.
\end{enumerate}
Hence, they ask for a simple way to determine whether the independence complex of $F$ is 
contractible and, if not, to compute the dimension of the associated sphere, and similarly to 
compute the dimension of the sphere associated to the dominance complex of $F$.

In this work we use the cores to answer these questions obtaining results which 
give topological interpretations to some well-known invariants of the underlying 
forest $F$ (e.g.~the domination number, the independent domination number, the matching number and the vertex covering number).
In particular, we prove that the 
contractibility of the independence complex of $F$ is detected by some graph theoretical properties of $F$. 
When the independence complex of $F$ is contractible, it collapses onto a point; 
when the independence complex of $F$ is not contractible, 
it collapses onto the boundary of a cross-polytope whose dimension equals 
the domination number of $F$.  Finally, we prove that the dominance complex of $F$ always 
collapses onto the boundary of a cross-polytope whose dimension equals the matching 
number of $F$.

The paper is organized as follows.

Section~\ref{sedue} contains the notation and background needed in the sequel. 
In Section~\ref{setre} we define the notion of domination between variables and we 
study its relationship to suspension.  
In Section~\ref{setreme} we introduce and study the main new concept of this work, namely 
the core of a simplicial complex; we reduce the study of the topology of a simplicial complex 
to the study of its core.  
In Section~\ref{sequa} we prove that the Euler characteristic of a simplicial complex can be easily 
computed from a set of generators of its Stanley-Reisner ideal. 
In Section~\ref{secin} we apply the method of the core to the independence complex $\Delta$ of 
a forest $F$. We find several conditions which are equivalent to the contractibility of $\Delta$, 
and we prove that, if  $\Delta$ is not contractible, then it collapses onto the boundary of a 
cross-polytope whose dimension equals the domination number of $F$.
In Section~\ref{sesei} we consider the problem of the independence complex of a general graph $G$. 
In Section~\ref{seset} we apply the method of the core to the dominance complex $\Delta$ of a forest 
$F$. We prove that $\Delta$ always collapses onto the boundary of a cross-polytope whose 
dimension equals the matching number of $F$.

\section{Notation and background} \label{sedue}

If $r \in \mathbb{Z}$, $r \geq 0$, we let $[r] := \{ 1 , \ldots , r\}$.  The cardinality of a set $A$ 
will be denoted by $|A|$.

We consider finite undirected graphs $G=(V,E)$ with no loops or multiple edges.  
For all $S \subset V$, let 
$N[S] := \bigl\{ w \in V \; | \; \exists s \in S , \{ s , w \} \in E \bigr\} \cup S$ be the 
{\it closed neighborhood of $S$}; when $S = \{ v \}$, then we let $N[v] = N[\{v\}]$.  A set 
$D \subset V$ is called {\it dominating} if for all $v \in V$, $N[v] \cap D \neq \emptyset $.  
A set $D \subset V$ is called {\it independent} if no two vertices in $S$ are adjacent, i.e.~$\{v,v'\}\notin E$ for all $v,v'\in D$.    
A {\it vertex cover of $G$} is a subset $C \subset V$ such that every edge of $G$ 
contains a vertex of $C$.  An {\it edge cover of $G$} is a subset $S \subset E$ such 
that the union of all the endpoints of the edges in $S$ is $V$.  
A {\it matching of $G$} is a subset $M \subset E$ of pairwise 
disjoint edges.

We consider the following classical invariants of a graph $G$ which have been extensively studied by graph theorists (see, for instance, 
\cite{AL}, \cite{ALH}, \cite{BC}, \cite{ET}, \cite{HHS}, \cite{HY}); we let 
\begin{itemize}
\item $\gamma (G) := \min \bigl\{ |D|, D \text{ is a dominating set of $G$} \bigr\}$ be the 
{\it domination number of $G$}; 
\item $i (G) := \min \bigl\{ |D|, D \text{ is an independent dominating set of $G$} \bigr\}$ 
be the {\it independent domination number of $G$}; 
\item $\alpha _0 (G) := \min \bigl\{ |C|, C \text{ is a vertex cover of $G$} \bigr\}$ be the 
{\it vertex covering number of $G$}; 
\item $\alpha _1 (G) := \min \bigl\{ |C|, C \text{ is an edge cover of $G$} \bigr\}$ be the 
{\it edge covering number of $G$}; 
\item $\beta _1 (G) := \max \bigl\{ |M|, M \text{ is a matching of $G$} \bigr\}$ be the 
{\it matching number of $G$}. 
\end{itemize}

Recall the following well-known results of K\"onig (cf~\cite{D}, Theorem~2.1.1) and Gallai 
(cf~\cite{HHS}, Theorem~9.27).

\begin{theorem}[K\"onig] \label{konig}
Let $G$ be a bipartite graph.  Then $\alpha _0 (G) = \beta _1 (G)$.
\end{theorem}

\begin{theorem}[Gallai] \label{gallai}
Let $G = (V,E)$ be a graph without isolated vertices.  Then 
$$ \alpha _1 (G) + \beta _1 (G) = |V| . $$
\end{theorem}

We refer the reader to~\cite{B} or~\cite{D} for all undefined notation on graph theory.

We let $X := \{ \ics \}$ and $\mathbb{Z}[X]$ be the polynomial ring with variables $\ics$ 
over the integers; we set $\mathbb{Z} [\emptyset ] := \mathbb{Z}$.  Let $m,m' \in \Rn$; 
we write $m' | m$ if $m'$ divides $m$.

\begin{definition}
A simplicial complex $\Delta $ on $X$ is a set of subsets of $X$, called {\it faces}, 
such that, if $\sigma \in \Delta $ and $\sigma ' \subset \sigma$, then $\sigma ' \in \Delta $.  
The faces of cardinality one are called {\it vertices}.

Equivalently, a simplicial complex $\Delta $ on $\Rn$ is a finite set of square-free monomials 
of $\Rn$ such that, if $m \in \Delta $ and $m'|m$, then $m' \in \Delta $.
\end{definition}

We do not require that $x \in \Delta $ for all $x \in X$.  We will frequently identify a set $S\subset X$ with the monomial $x_1^{\varepsilon_1}\cdots x_n^{\varepsilon_n}$, where $\varepsilon_i=\left\{ 
\begin{array}{ll}
1, & \textrm{if $x_i\in S$;}\\
0, & \textrm{if $x_i\notin S$.}
\end{array} \right.
$
Note that the empty set is identified with the monomial $1$. 
We refer the reader 
to~\cite{MS} for all undefined concepts from commutative algebra.

Every simplicial complex $\Delta $ on $\Rn$ different from $\{1\}$ has a standard geometric 
realization.  Let $e_1 , \ldots , e_n$ 
be the standard basis of $\mathbb{R}^n$.  The realization of $\Delta $ is the union of the convex 
hulls of the sets $\{e_i {\text{ such that }} x_i | m \}$, for each monomial $m \in \Delta $.  Whenever 
we mention a topological property of $\Delta $, we implicitly refer to the geometric realization 
of $\Delta $.

Let $I \subset \Rn$ be a monomial ideal (i.e.~an ideal generated by monomials) containing 
$\icss$.  The set of monomials 
of $\Rn \setminus I$ is a simplicial complex on $\Rn$ that we denote by $R(I)$.  Conversely, given 
a simplicial complex $\Delta $ on $\Rn$, let $I _\Delta \subset \Rn$ be the ideal 
generated by the monomials not in $\Delta $.  Clearly $\Delta = R (I_\Delta )$ and 
$I = I_{R(I)}$.  Note that $I_\Delta $ is (essentially) the Stanley-Reisner ideal of the 
simplicial complex $\Delta $ (see~\cite{S}).

As examples, consider the ideals $I_n = (\icss )$, $J_n = (x_1 \cdots x_n , \icss )$, 
and $K_n = (x_1 x_2 , x_3 x_4 , \ldots , x_{2n-1} x_{2n} , x_1^2 , \ldots , x_{2n} ^2 )$.  
The simplicial complex $R(I_n)$ is the $(n-1)-$dimensional simplex and $R(J_n)$ is 
its boundary; $R(K_n)$ is the boundary of the $n-$dimensional cross-polytope, 
which is the dual of the $n-$dimensional cube.  Note that the cube, its boundary and 
the cross-polytope are not simplicial complexes.  Furthermore $R(K_n)$ is the 
$n-$th suspension of the simplicial complex $\{ 1 \}$.

From now on, unless explicitly mentioned otherwise, $I$ denotes a monomial ideal 
of $\Rn$ containing $\icss $, and $\Delta := R(I)$.

For the basic concepts of simplicial homology we refer the reader to~\cite{Mu}.  
We identify the free abelian group of simplicial chains on $\Delta $ with the quotient $\Rn / I$: 
the chains of dimension $i-1$ are the span of the monomials of degree $i$.  Choose 
an order on $X$; this induces a boundary map $\delta $ on the simplicial chains.  
We denote by $Z(\Delta)$ the $\mathbb{Z}-$module of cycles on $\Delta $ and by 
${\rm \tilde H} \bigl( \Delta , \mathbb{Z} \bigr)$ the reduced homology groups with integer 
coefficients of $\Delta $.  A {\it quasi-isomorphism of degree $r$} is a morphism of chains 
sending chains of degree $k$ to chains of degree $k+r$ which induces an isomorphism 
in homology.

Let $M$ be a finitely generated graded $\Rn-$module $M$; we denote  the (multi-graded) Hilbert series of $M$ by $H \bigl( M ; \ics \bigr)$. The {\it multi-graded face polynomial $\cF_{\Delta } (\ics)$ of 
$\Delta $} is the polynomial of $\Rn$ 
$$ \cF_{\Delta } (\ics) := \sum _{m \in \Delta } m = 
H \bigl( \Rn/ I ; \ics \bigr) . $$
The {\it face polynomial $F_\Delta (t)$ of $\Delta $} is the polynomial 
$$ F_\Delta (t) := \sum _{m \in \Delta } t^{\deg m} = \cF_\Delta (t , \ldots , t) . $$
The {\it reduced Euler characteristic of $\Delta $} is 
$\tilde e (\Delta ) := - F_\Delta (-1)$.

We note that the simplicial complexes $R \bigl( (\ics) \bigr) = \{1\}$ and 
$R \bigl( \Rn \bigr) = \emptyset $ are different: we call $\{ 1 \}$ the $(-1)-$dimensional 
sphere, and $\emptyset $ the $(-1)-$dimensional simplex.  The empty simplex 
$R \bigl( \Rn \bigr)$ is contractible.  For $n \geq 1$, let 
$S^{n-2} := R \bigl( (x_1 \cdots x_n , \icss ) \bigr)$, 
the {\it sphere of dimension $n-2$}.  Consistently with these conventions, the reduced 
Euler characteristic of the $(-1)-$dimensional sphere is $-1$ while the reduced Euler 
characteristic of the $(-1)-$dimensional simplex is $0$.

Let $x \in \Rn$ be a monomial and define simplicial complexes 
$$ \begin{array}{cccccc}
(\Delta : x) & := & \bigl\{ m \in \Delta \; | \; \; xm \in \Delta \bigr\} & = & 
R \bigl( I : x \bigr) \\[5pt]
(\Delta , x) & := & \bigl\{ m \in \Delta \; | \; \; x \nmid m \bigr\} & = & 
R \bigl( I , x \bigr) , 
\end{array} $$
where $(I:x) = \{ m \in \Rn \; | \; xm \in I \}$ and $(I,x)$ is the ideal generated by $I$ and $x$. The simplicial complexes $(\Delta : x)$ and $(\Delta , x)$ are usually called  link and face-deletion of $x$.
If $I _1 , \ldots , I_k \subset \Rn$ are monomial ideals containing $\icss$, then we define 
$$ {\rm join} \bigl( R (I_1) , \ldots , R (I_k) \bigr) := 
\bigl\{ \lcm \{m_i , i \in [k] \} \; | \; m_i \in R (I_i) \bigr\} .  $$
If $x$ and $y$ are monomials, let 
$$ \begin{array}{rcl}
A_x \bigl( \Delta \bigr) & := & {\rm join} \bigl(\Delta, \{1,x\} \bigr) \\[5pt]
\Sigma_{x,y} \bigl( \Delta \bigr) & := & 
{\rm join} \bigl(\Delta, \{1,x,y\} \bigr) . 
\end{array} $$

If $x,y \in X$
then 
$A_x \bigl( \Delta \bigr)$ and $\Sigma_{x,y} \bigl( \Delta \bigr)$ 
are both simplicial 
complexes.
If $x \neq y \in X$  and they are coprime with the faces of $\Delta $, then 
$A_x \bigl( \Delta \bigr)$ and $\Sigma_{x,y} \bigl( \Delta \bigr)$ are called respectively the {\it cone on $\Delta$ with apex $x$} and the {\it suspension of $\Delta$}.  If $x \neq y$ and $x' \neq y'$ are 
variables in $X$ coprime with all the faces of $\Delta $, then the suspensions 
$\Sigma_{x,y} \bigl( \Delta \bigr)$ and $\Sigma_{x',y'} \bigl( \Delta \bigr)$ are 
isomorphic; hence in this case sometimes we drop the subscript 
from the notation.
It is well-known that 
if $\Delta $ is contractible, then $\Sigma (\Delta )$ is contractible, and that 
if $\Delta $ is homotopic to a sphere of dimension $k$, then $\Sigma (\Delta )$ is 
homotopic to a sphere of dimension $k+1$.

We recall the notions of collapse and simple-homotopy (see \cite{C}).  
Let $\sigma \supset \tau $ be faces of a simplicial complex 
$\Delta $ and suppose that $\sigma $ is maximal and $\deg (\tau ) = \deg (\sigma ) -1$.  
If $\sigma $ is the only face of $\Delta $ properly containing $\tau $, then the removal 
of $\sigma $ and $\tau $ is called an {\it elementary collapse}. If a simplicial complex $\Delta '$ is obtained from $\Delta$ by an elementary collapse,  we write 
$\Delta \succ \Delta'$.

Equivalently in terms of ideals, an elementary collapse is obtained by adding to the monomial 
ideal $I$ a monomial $\tau $ such that 
\begin{itemize}
\item $\tau $ is a monomial not in $I$, and 
\item there is a unique variable $a$ such that $\sigma := a \tau $ is also not in $I$.
\end{itemize}

When $\Delta '$ is a subcomplex of $\Delta $, we say that {\it $\Delta $ collapses onto $\Delta '$} 
if there is a sequence of elementary collapses leading from $\Delta $ to $\Delta '$.

\begin{definition}
Two simplicial complexes $\Delta $ and $\Delta '$ are {\it simple-homotopic} if they are 
equivalent under the equivalence relation generated by $\succ $.
\end{definition}

It is clear that if $\Delta $ and $\Delta '$ are simple-homotopic, then they are also homotopic,
and that a cone collapses onto a point.

Let $a \in X$.  If $a \in I$, then $a$ is not a vertex of $\Delta $.  Since we are 
interested in studying $\Delta $, we identify $I \subset \Rn$ with 
$J \subset \mathbb{Z} [X \cup \{ a \}]$ whenever $J  = (I,a)$, because the associated 
simplicial complexes on $X$ and $X \cup \{ a \}$ are the same.  Note that, in general, 
a monomial ideal $J$ has a unique minimal generating set $M$ consisting of 
monomials.  If $J = I_\Delta $, then we let $M = B \cup \{ \icss \}$ with 
$B \cap \{ \icss \} = \emptyset $ and we call the elements 
of $B$ the minimal square-free generators of $I _\Delta $.  It follows from the definitions 
that $(I : a) = (I,a)$ if and only if $\Delta $ is a cone with apex $a$; equivalently, $\Delta $ 
is a cone with apex $a$ if and only if $a$ divides no monomial of $B$.

\begin{lemma} \label{precontra}
Let $\Delta $ be a simplicial complex and let $x \in \Rn$ be a monomial; then 
$\Delta = A_x (\Delta : x) \cup (\Delta , x)$.
\end{lemma}

\begin{proof}
There is an exact sequence 
$$ 0 \to \Rn/(I : x) \stackrel {\cdot x} {\longrightarrow} \Rn/ I \longrightarrow 
\Rn/( I , x) \to 0 $$
and hence $\cF _{\Delta} = x \cF _{\left( \Delta : x \right)} + \cF _{\left( \Delta , x \right)}$.  
On the other hand the multi-graded face polynomial of 
$A_x (\Delta : x) \bigcup _{(\Delta : x)} (\Delta , x)$ is 
$$ \cF _{(\Delta : x)} + x \cF _{(\Delta : x)} + \cF _{(\Delta , x)} - \cF _{(\Delta : x)} . $$
Note that $x$ is coprime with all vertices of $(\Delta : x)$.
\end{proof}

\section{Domination} \label{setre}

In this section we introduce the notion of domination between variables in a monomial 
ideal $I \subset \Rn $ containing $\icss$, and we give some preliminary results on the 
topology of the simplicial complex $R(I)$.

Since domination and suspension are closely related, we start with some remarks on 
suspensions of simplicial complexes.  
It is immediate that the suspension $\Sigma_{x,y} \Delta$ is a cone with apex $a$ 
if the simplicial complex $\Delta$ is a cone with apex $a$, independently of whether 
or not $x$ and $y$ 
are vertices of $\Delta$. On the contrary, if two simplicial complexes are 
simple-homotopic, it does not follow in general that their suspensions are 
homotopic. The next lemma analyses the question of lifting collapses to suspensions.

\begin{lemma} \label{casetti}
Let $\Delta \succ \Delta '$ be simplicial complexes, $\sigma \supset \tau $ the faces 
removed in the elementary collapse, and $x ,y \in X$, $x \neq y$.  Then 
the simplicial complex $\Sigma _{x,y} \bigl( \Delta )$ collapses onto 
$\Sigma _{x,y} \bigl( \Delta ' )$ unless one of the following is satisfied: 
\begin{enumerate}
\item $x | \tau $, $y \nmid \sigma $ and $y \frac{\sigma } {x} \notin \Delta '$, $y \frac{\tau } {x} \in \Delta '$; 
\label{cauno}

\smallskip

\item $y | \tau $, $x \nmid \sigma $ and $x \frac{\sigma } {y} \notin \Delta '$, $x \frac{\tau } {y} \in \Delta '$.
\label{cadue}
\end{enumerate}
In these last cases the two suspensions have different Euler characteristics.
\end{lemma}

\begin{proof}
For notational convenience we may write $a$ for the singleton $\{ a \}$.

Notice that $\Delta = \Delta ' \cup \tau \cup \sigma $ and hence 
\begin{equation} \label{unioni}
\Sigma _{x,y} \bigl( \Delta ) = \Sigma _{x,y} \bigl( \Delta ' ) \cup \tau \cup \sigma 
\cup \lcm \{ x , \tau \} \cup \lcm \{ x , \sigma \} 
\cup \lcm \{ y , \tau \} \cup \lcm \{ y , \sigma \}
\end{equation}
%
%
where the unions need not be disjoint.  We separate six mutually exclusive cases.

{\bf Case 1.}  $x,y | \sigma $.  We show that 
$\Sigma _{x,y} \bigl( \Delta ) = \Sigma _{x,y} \bigl( \Delta ' )$.  By (\ref{unioni}) it 
suffices to check that $\sigma \in \Sigma _{x,y} \bigl( \Delta ' )$.  By hypothesis 
there is $t \in \{ x , y \}$ such that $t | \tau $; then $\frac {\sigma } {t} \in \Delta '$ and 
hence $\sigma \in \Sigma _{x,y} \bigl( \Delta ' )$.

{\bf Case 2.}  $x,y \nmid \sigma $.  The union (\ref{unioni}) is disjoint and 
$\Sigma _{x,y} \bigl( \Delta ' )$ is obtained from $\Sigma _{x,y} \bigl( \Delta )$ by the 
elementary collapses of the faces $x \sigma \supset x \tau $ and $y \sigma \supset y \tau $ 
(in any order), followed by the collapse of the faces $\sigma \supset \tau $.

{\bf Case 3.}  $\sigma = x \tau $, $y \nmid \sigma $.  By (\ref{unioni}) we have 
$\Sigma _{x,y} \bigl( \Delta ) = 
\Sigma _{x,y} \bigl( \Delta ' ) \cup \tau \cup \sigma \cup y \tau \cup y \sigma $; collapsing 
successively the faces $y \sigma \supset y \tau $ and $\sigma \supset \tau $ we conclude.

{\bf Case 4.}  $\sigma = y \tau $, $x \nmid \sigma $.  Follows by symmetry from Case 3.

{\bf Case 5.}  $x | \tau $, $y \nmid \sigma $.  Note that $\frac {\sigma } {x} \in \Delta '$ and 
hence $\sigma \in \Sigma _{x,y} \bigl( \Delta ' )$.  Thus by (\ref{unioni}) we have 
$\Sigma _{x,y} \bigl( \Delta ) = 
\Sigma _{x,y} \bigl( \Delta ' ) \cup y \tau \cup y \sigma $.  We have 
$$ \begin{array} {ccc} 
y \sigma  \in \Sigma _{x,y} \bigl( \Delta ' ) & \Longleftrightarrow & y \frac {\sigma } {x} \in \Delta ' \\[5pt]
\Downarrow & & \Downarrow \\[5pt]
y \tau \in \Sigma _{x,y} \bigl( \Delta ' ) & \Longleftrightarrow & y \frac {\tau } {x} \in \Delta ' 
\end{array} $$
and hence 
\begin{itemize}
\item if $y \frac {\sigma } {x} \in \Delta '$, then $\Sigma _{x,y} \bigl( \Delta ) = 
\Sigma _{x,y} \bigl( \Delta ' )$; 

\smallskip

\item if $y \frac {\tau } {x} \notin \Delta '$, then $\Sigma _{x,y} \bigl( \Delta ' )$ is obtained 
by the elementary collapse of the faces 
$y \sigma \supset y \tau $ of $\Sigma _{x,y} \bigl( \Delta )$; 

\smallskip

\item if $y \frac {\sigma } {x} \notin \Delta '$ and $y \frac {\tau } {x} \in \Delta '$, then (\ref{cauno}) is 
satisfied.  Note that $\Sigma _{x,y} \bigl( \Delta ) = \Sigma _{x,y} \bigl( \Delta ' ) \cup y \sigma $; 
thus the Euler characteristics of $\Sigma _{x,y} \bigl( \Delta )$ and $\Sigma _{x,y} \bigl( \Delta ' )$ 
differ by one.
\end{itemize}

{\bf Case 6.}  $y | \tau $, $x \nmid \sigma $.  Follows by symmetry from Case 5.
\end{proof}

\begin{remark}\label{casotti}
Note that if at least one among $x$ or $y$ is not a vertex of $\Delta '$, then  $\Sigma _{x,y} \bigl( \Delta )$ collapses onto 
$\Sigma _{x,y} \bigl( \Delta ' )$ since (\ref{cauno}) and (\ref{cadue}) cannot be satisfied.
\end{remark}

\begin{example}
Consider the following ideals of $\mathbb{Z} [x,y,u,v]$: 
\begin{itemize}
\item $I = (xy , yu , x^2 , y^2 , u^2 , v^2)$, 
\item $J = (xy , yu , xv , x^2 , y^2 , u^2 , v^2)$.
\end{itemize}
\begin{figure}[ht]
\psfrag{x}{\small $x$}
\psfrag{y}{\small $y$}
\psfrag{v}{\small $v$}
\psfrag{I}{\small $R(I)$}
\psfrag{J}{\small $R(J)$}
\psfrag{u}{\small $u$}
\begin{center}
{\includegraphics[scale=0.5]{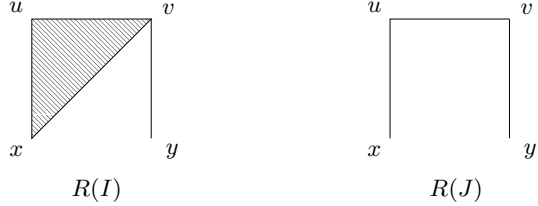}} \par
\caption{\label{trim2.fig} The simplicial complexes $R(I)$ and $R(J)$}
\end{center}
\end{figure}

The complex $R(I)$ collapses onto $R(J)$ by the elementary collapse of the faces 
$xuv \supset xv$.  On the other hand $\Sigma _{x,y} \bigl( R(I) \bigr)$ is the three 
dimensional simplex with vertices $x,y,u,v$, while $\Sigma _{x,y} \bigl( R(J) \bigr)$ 
is its boundary.  This is case~(\ref{cauno}) of Lemma~\ref{casetti}.
\end{example}

We now give the main definition of this section.

\begin{definition}
Let $a,b \in X$; {\it $a$ dominates $b$ in $I$} 
if $R(I)$ is not a cone with apex $b$ and $R(I,a)$ is a cone with apex $b$.  
\end{definition}

\begin{figure}[ht]
\psfrag{a}{\small $a$}
\psfrag{b}{\small $b$}
\psfrag{punti}{\small $R(I,a)$}
\psfrag{virgola}{\small $R(I:a)$}
\begin{center}
{\includegraphics[scale=0.5]{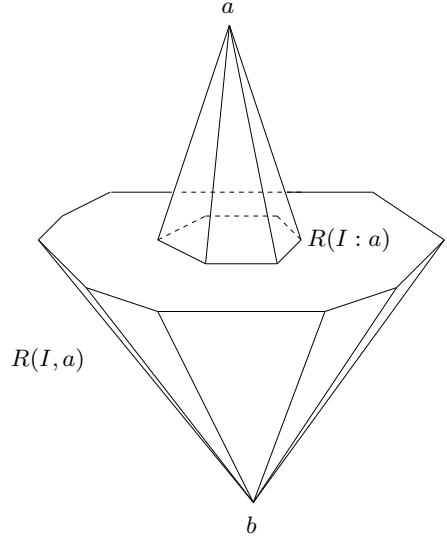}} \par
\caption{\label{trim.fig} The variable $a$ dominates $b$}
\end{center}
\end{figure}

Note that $a$ dominates $b$ in $I$ if and only if every minimal square-free generator of $I$ 
divisible by $b$ is also divisible by $a$ and there are such monomials. 
Loosely speaking, if $a$ dominates $b$ then $R(I)$ is composed out of two cones:
the cone $R(I,a)$ with apex $b$ and the cone with apex $a$ on the subcomplex  $R(I:a)\subset R(I,a)$ (see Figure~\ref{trim.fig}).
The apex $a$ is not a vertex of $R(I,a)$ while $b$ might be a vertex of $R(I:a)$.

An immediate consequence of Lemma~\ref{precontra} is that, 
if $a,b \in X$ and $a$ dominates $b$ in $I$, then $R(I)$ is 
homotopic to $\Sigma \bigl( R(I : a) \bigr)$.  In fact we can prove that $R(I)$ is 
simple-homotopic to $\Sigma \bigl( R(I : a) \bigr)$.


\begin{theorem} \label{paba}
Let $a$ dominate $b$ in $I$ and let $I' = \bigl( I, a b', (b')^2 \bigr) \subset \Rn [b']$.  Then 
\begin{enumerate}
\item $R(I')$ collapses onto $R(I)$, \label{semplomo2}
\item $R(I)$ collapses onto $\Sigma _{a,b} \bigl( R(I:a) \bigr)$, \label{paba2}
\item $R(I')$ collapses onto $\Sigma _{a,b'} \bigl( R(I':a) \bigr)$. \label{paba3}
\end{enumerate}
In particular $R(I)$ is simple-homotopic to $\Sigma \bigl( R(I:a) \bigr)$.
$$ \xymatrix{
& R(I') \ar[dl] _{collapses} \ar[dr] ^{collapses} \\
R(I) \ar[d] _{collapses} & & \Sigma _{a,b'} \bigl( R(I':a) \bigr) \ar@{=} [d] \\
\Sigma _{a,b} \bigl( R(I:a) \bigr) & & \Sigma \bigl( R(I:a) \bigr)
} $$
\end{theorem}

\begin{proof}
(\ref{semplomo2}) By definition $R(I) \subset R(I')$.  We show that there is a sequence 
$\Delta _0 := R(I') \succ \Delta _1 \succ \cdots \succ \Delta _s := R(I)$ of simplicial complexes 
where $s$ is the number of faces of $R(I')$ containing $b$ and $b'$, and $\Delta _i$ is 
obtained from $\Delta _{i-1}$ by the elementary collapse of a face $\sigma $ containing $b$ 
and $b'$ and the face $\tau = \frac {\sigma } {b}$.  Let $\Delta _0 := R(I')$; suppose that 
$\Delta _j$ has been defined for all $j \leq i \leq s$.

If $i=s$, then $b b'$ is not a face of $\Delta _i$, since at each step we remove 
exactly one face containing $b b'$.  In this case we are done, since we have already removed 
the faces $\sigma = bb'$ and $\frac {\sigma } {b} = b'$, and thus $\Delta _s = R(I)$.

If $i < s$, we define $\Delta _{i+1}$ as follows.  Note that $b b'$ is a face of $\Delta _i$, since 
we removed fewer than $s$ faces containing $b b'$, and let $f bb'$ be a maximal face of 
$\Delta _i$ containing $bb'$.  We prove first that the only face strictly containing $f b'$ is $f bb'$.  
Let $q f b'$ be a maximal face of $\Delta _i$ containing $f b'$ and assume by contradiction 
that $q f b' \neq f bb'$.  Clearly $a$ does not divide $q f$ since $a b' \in I'$, and $b$ cannot 
divide $q$ since $f bb'$ is maximal.  By the assumption 
on the elementary collapses, the faces of $R(I')$ not containing $b'$ are not affected by the collapses.  
Hence $q f b$ is a face of $R(I,a)$, since $a$ dominates $b$ and $a \nmid qf$, and it is 
also a face of $\Delta _j$ for all $j \leq i$; the monomial $q f bb'$ is a 
face of $\Delta _0$ and not of $\Delta _i$ by maximality.  Hence there is an index $j<i$ such that 
$\Delta _{j+1}$ is obtained by removing the faces $\sigma = q f bb'$ and $\frac {\sigma } {b}$, 
contradicting the fact that $q f b' \in \Delta _i$.

By what we just proved, we may collapse the faces $q f bb'$ and $q f b'$.  We define 
$\Delta _{i+1}$ to be the result of this collapse.  Iterating this procedure we conclude.

\noindent
(\ref{paba2})  The simplicial complex $\Sigma _{a,b} \bigl( R (I:a) \bigr)$ is a subcomplex of $R(I)$.  
Let $\sigma $ be a maximal face of $R(I)$ not in $\Sigma _{a,b} \bigl( R (I:a) \bigr)$; we show 
that $\sigma $ contains $b$ and that we may collapse $\sigma $ and $\frac {\sigma } {b}$.  
Note that $\sigma $ does not contain $a$, since $\sigma $ is not a face of 
$\Sigma _{a,b} \bigl( R (I:a) \bigr)$, and hence 
$\sigma $ is a face of $R(I,a)$; since $R(I,a)$ is a cone of apex $b$ and $\sigma $ is maximal, 
$\sigma $ contains $b$.  Write $\sigma  = \tau b$; if $\tau a$ is a face of $R(I)$, then 
$\tau \in R(I:a)$ and hence $\tau b \in \Sigma _{a,b} \bigl( R (I:a) \bigr)$.  Since this is not the 
case, if $\tau h$ is a face of $R(I)$ containing $\tau $, then $h$ is not divisible by $a$ 
and hence $\tau h$ is a face in $R(I,a)$; thus $\lcm \{ \tau h , b \}$ is also a face of $R(I,a)$ and by 
maximality of $\sigma = \tau b$ we conclude that $h | b$.  Thus the only face of $R(I)$ strictly 
containing $\tau $ is $\sigma $.  Hence we may collapse the faces $\sigma $ and $\tau $ to obtain 
a simplicial complex $\Delta '$.  Note that $\bigl( \Delta ' , a \bigr)$ is again a cone with apex $b$ 
and we may iterate this procedure and conclude.

\noindent
(\ref{paba3})  Follows from part~(\ref{paba2}), since $a$ dominates $b'$ in $I'$.

The last statement follows since $\Sigma _{a,b'} \bigl( R(I':a) \bigr)$ is isomorphic 
to the abstract suspension $\Sigma \bigl( R(I:a) \bigr)$ because $a,b' \notin R(I':a) = R(I:a)$.
\end{proof}

\begin{example}
Let $I = (x_1 x_2 x_3 , x_1^2 , x_2^2 , x_3^2)$ and 
$J = (x_1 x_2 , x_3 x_4 , x_1^2 , x_2^2 , x_3^2 , x_4^2)$.  We have that $x_3$ dominates 
$x_1$ in $I$ and $x_3$ dominates $x_4$ in $J$; moreover 
$R(I:x_3) = R(J:x_3)$ is the simplicial complex consisting of the two points $x_1$, $x_2$.  
Hence both $R(I)$ and $R(J)$ are simple-homotopic to the boundary of the $2-$dimensional 
cross-polytope (a square) by Theorem~\ref{paba}.  Whereas the simplicial complex $R(J)$ 
is actually the boundary of the $2-$dimensional cross-polytope, the simplicial complex $R(I)$ is the 
boundary of the $2-$dimensional simplex (a triangle).
\end{example}

The following lemmas are needed in the next section.

\begin{lemma} \label{elcontre}
Let $\Delta \succ \Delta '$ be simplicial complexes and let $\sigma \supset \tau $ be the faces 
removed in the elementary collapse.  Then 
\begin{enumerate}
\item any cycle of $\Delta $ is a combination of faces different from $\sigma $; \label{sig}
\item any cycle of $\Delta $ is homologous to a combination of faces different from $\tau $; \label{omota}
\item the inclusion $\Delta ' \subset \Delta $ induces an isomorphism in homology. \label{inclu}
\end{enumerate}
\end{lemma}

\begin{proof}
Write $\sigma = a \tau $ and choose an order of the variables such that $a$ is first.

\noindent
(\ref{sig}) Let $z = c \sigma + \sum _{A \neq \sigma } c_A A$ be a cycle.  Since 
$$ 0 = \partial z = c \tau - c a \partial \tau + \sum _{A \neq \sigma } c_A \partial A $$
we have $c=0$, because the face $\tau $ is properly contained only in $\sigma $.

\noindent
(\ref{omota}) Let $z = d \tau + \sum _{A \neq \tau } c_A A$ and note that $z - \partial (d \sigma )$ 
has the required property.

\noindent
(\ref{inclu}) Let $z = \partial \bigl( c \sigma + d \tau + \sum _{A \neq \sigma , \tau } c_A A \bigr)$.  
By part (\ref{omota}) we may assume that the coefficient of $\tau $ in $z$ is zero and hence 
that $c=0$.  Then $z$ is the boundary of 
$-d \partial \sigma + d \tau + \sum _{A \neq \sigma , \tau } c_A A$ and we are done.
\end{proof}

\begin{remark} \label{cretino}
If $\Delta '$ is obtained from $\Delta $ by a sequence of elementary collapses, 
then the inclusion $\Delta ' \subset \Delta $ induces an isomorphism in homology.
\end{remark}

\begin{lemma} \label{separati}
Let $a$ dominate $b$ in $I$ and let $I' = \bigl( I, a b', (b')^2 \bigr) \subset \Rn [b']$.  Then 
the inclusion $R(I) \subset R(I')$ induces an isomorphism in homology whose inverse 
is induced by the map $\bar \varphi $ of chains given by 
\smallskip
$$ \begin{array}{rcl}
m & \longmapsto & \left\{ \begin{array}{ll}
m, & {\text{ if $m$ does not contain $b'$,}} \\[3pt]
b \frac{m}{b'}, & {\text{ if $m$ contains $b'$ and does not contain $b$,}} \\[3pt]
0, & {\text{ if $m$ contains $b b'$,}}
\end{array}
\right.
\end{array} $$
where $m$ is any face of $R(I')$, with the variables ordered so that $a<b'<b<x$ for all 
remaining variables $x$.
\end{lemma}

\begin{proof}
First of all we check that $\bar \varphi $ is a map of chains.  Let $m \in R(I')$ 
be a face; we only need to consider the case $m=b' m'$, with $m'$ not containing $b$.  
In this case, $m'$ cannot contain $a$ since $a b'$ is not a face of $R(I')$, and, because 
$a$ dominates $b$ in $I$, it follows that $b m'$ is a face of $R(I)$, as needed.  Let $\bar \iota $ be the 
map of chains induced by the inclusion $R(I) \subset R(I')$ and let $z$ be a cycle in $R(I')$.  
To conclude it is enough to check that $z - \bar \iota \bigl( \bar \varphi (z) \bigr)$ is a boundary 
in $R(I')$.  We may write $z = b' b A + b' B + C$, where the chains in $A,B$ do not contain $b,b'$ 
and the chains in $C$ do not contain $b'$; note that the chains in $A,B$ cannot contain $a$ 
since $a b'$ is not a face.  Since $z$ is a cycle we have 
$$ 0 = \partial z = bA - b'A + b'b\partial A + B - b' \partial B + \partial C $$
which implies that $-(A+\partial B) = 0$ since it is the coefficient of the faces containing $b'$ 
and not containing $b$.  Hence we may write $z - \bar \iota \bigl( \bar \varphi (z) \bigr) = 
b'bA + (b'-b)B = \partial (b'bB)$.  Note that $b'bB$ is a chain since the faces in $B$ do not 
contain $a$ and $bB$ is a chain.
\end{proof}

We note that applying $\bar \varphi $ to a chain simply deletes all terms containing the 
face $bb'$ and replaces $b'$ by $b$ in all remaining terms.

\section{Resolutions and cores} \label{setreme}

In this section we introduce the notions of resolution and core of a monomial 
ideal $I \subset \Rn $ containing $\icss$, and we deduce topological properties 
of the simplicial complex $R(I)$ from the resolution and the core of $I$.

Let $( a_1 , \ldots , a_r )$ be a sequence of variables of $\Rn$ and, for $i \in [r+1]$, 
let $I_i := (I : a_1 \cdots a_{i-1})$.  Consider the following properties: 
\begin{enumerate}
\item  for $i \in [r]$ we have $a_i \notin I_i$; \label{inizio}
\item for all $i \in [r]$ we have either $R(I_i)$ is a cone with apex $a_i$ 
or there exists $b_i \in X$ such that $a_i$ dominates $b_i$ in $I_i$. \label{sec}
\end{enumerate}

\begin{definition}
A {\it resolution of $I$} is a sequence $A=( a_1 , \ldots , a_r )$ satisfying (\ref{inizio}) 
and (\ref{sec}).  We call $c(A) := I_{r+1}$ the {\it core of $A$}, $d(A):=r$ the 
{\it depth of $A$}, and 
$$\begin{array}{rcl}
c(I)&:=&\{c(A)|\text{$A$ is maximal}\}\\[5pt]
d(I)&:=&\min \{d(A)|\text{$A$ is maximal}\}
\end{array}$$
respectively the {\it core} and the {\it depth of $I$}.  

The resolution $A$ is {\it spherical} if the simplicial complex $R(I_i)$ is not a cone of 
apex $a_i$, for all $i \in [r+1]$. The ideal $I$ is {\it spherical} if it admits a maximal 
resolution which is spherical, {\it conical} if it admits a resolution which is not spherical, 
{\it simple} if $(\ics) \in c(I)$.
\end{definition}

\begin{remark} \label{sfeco}
We shall see that the properties of being spherical and conical are mutually exclusive 
(cf Theorem~\ref{coco}).
\end{remark}

Note that if $R(I)$ is a cone with apex $b$ and $a \neq b$, then $R(I:a)$ is a cone with 
apex $b$.  Hence if $A$ is a resolution of $I$ and $R \bigl( c(A) \bigr)$ is a cone with 
apex $b$, then all maximal resolutions extending $A$ must contain $b$.

\begin{theorem} \label{collacolla}
Let $A=(a_1 , \ldots , a_r)$ be a resolution of $I$.  
\begin{itemize}
\item If $A$ is conical, then $R(I)$ collapses onto a point.
\item If $A$ is spherical and $a_i$ dominates $b_i$ in $(I : a_1 \cdots a_{i-1})$, then 
$R(I)$ collapses onto ${\rm join} \bigl( \Sigma , R(c(A)) \bigr)$, where 
$\Sigma := {\rm join} \bigl( \{1 , a_1 , b_1\} , \ldots , \{1 , a_r , b_r\} \bigr)$.  In particular $R(I)$ is 
simple-homotopic to $\Sigma ^{d(A)} \bigl( R(c(A)) \bigr)$.
\end{itemize}
\end{theorem}

\begin{proof}
If $R(I)$ is a cone of apex $a_1$, then it collapses onto the point $a_1$.  Otherwise by 
Theorem~\ref{paba} the simplicial complex $R(I)$ collapses onto 
$\Sigma _{a_1,b_1} R \bigl( I : a_1 \bigr)$.  


Suppose first that $A$ is spherical and proceed by induction on $r$. 
By Theorem~\ref{paba} part~(\ref{paba2}), the simplicial complex $R(I)$ collapses onto 
$\Sigma _{a_1,b_1} R \bigl( I : a_1 \bigr)$.  If $r=1$, we are done. 
Suppose $r\geq 2$. 
Let $A'=(a_2 , \ldots , a_r)$ and $\Sigma' := {\rm join} \bigl( \{1 , a_2 , b_2\} , \ldots , \{1 , a_r , b_r\} \bigr)$.  Note that $c(A')=c(A)$.
To conclude, it suffices to show that $\Sigma _{a_1,b_1} R \bigl( I : a_1 \bigr)$ collapses onto $\Sigma _{a_1,b_1} {\rm join} \bigl( \Sigma' , R(c(A')) \bigr)={\rm join} \bigl( \Sigma , R(c(A)) \bigr)$.

By induction, $R \bigl( I : a_1 \bigr)$ collapses onto ${\rm join} \bigl( \Sigma' , R(c(A')) \bigr)$, since $A'$ is a spherical resolution of $(I:a_1)$. Since $a_1$ is not a vertex of $R(I:a_1)$, we apply repeatedly Lemma~\ref{casetti} (see Remark~\ref{casotti}) to conclude.

Suppose now that $A$ is conical. Let $i$ be the smallest index such that $R(I:a_1\cdots a_{i-1})$ is a cone with apex $a_i$. By what we just proved, $R(I)$ collapses onto $C := \Sigma_{a_1 , b_1} \cdots \Sigma_{a_{1-1} , b_{i-1}} R(I:a_1\cdots a_{i-1})$.  Since $C$ is an iterated suspension of a cone, it is a cone; thus $C$ and hence $R(I)$ collapse to a point.
\end{proof}



\begin{remark} \label{recolla}
Let $a_1 , \ldots , a_r$ be distinct variables and let $b_1 , \ldots , b_r$ be variables such 
that $a_j \neq b_i$ for all $j \leq i$.  Hence 
$$ r+1 \leq s := |\{ a_1 , \ldots , a_r , b_1 , \ldots , b_r\}| \leq 2r . $$
Let $\Sigma := \Sigma _{a_1,b_1} \cdots \Sigma _{a_r,b_r} (\{1\}) $.  
If $s=2r$, then $\Sigma $ is the 
boundary of the $r-$dimensional cross-polytope.  If $s = r+1$, then $\Sigma $ is 
the boundary of the $r-$dimensional simplex.  In the cases in which $r+1 < s < 2r$, 
the complex $\Sigma $ is a hybrid of the two extreme cases.  All the complexes 
thus obtained are simple-homotopic.



The first case in which $\Sigma $ may be different from the boundary of a simplex or 
a cross-polytope is for $r=3$.  Let $I=(x_1 x_2 , x_3 x_4 x_5 , x_1^2 , \ldots , x_5^2)$ and 
let $(a_1,a_2,a_3) = (x_1,x_3,x_4)$ and $(b_1,b_2,b_3) = (x_2,x_5,x_5)$.  We have 
$a_i$ dominates $b_i$ in $(I:a_1 \cdots a_{i-1})$ and $s=5$; the simplicial complex 
$R(I)$ is the suspension of the boundary of a triangle.
\end{remark}

The following result gives an explicit description of the homology of a simplicial complex 
$\Delta $ in terms of the homology of the core of a resolution of $\Delta $.

%

\begin{theorem} \label{omolomap}
Let $A = (a_1 , \ldots , a_r)$ be a spherical resolution of $I$ and suppose that 
$a_i$ dominates $b_i$ in $I_i := (I:a_1 \cdots a_{i-1})$.  Then there is a quasi-isomorphism 
$\varphi : Z \bigl( R(c(A)) \bigr) \rightarrow Z \bigl( R(I) \bigr)$ of degree $r$.

Furthermore $\varphi = \pi \circ \varphi ' $ where $\pi$ is a map sending each face 
$\sigma $ to $\pm \sigma $ and 
$$ \begin{array}{rcl}
\varphi ' : Z \bigl( R(c(A)) \bigr) & \longrightarrow & Z \bigl( R(I) \bigr) \\[5pt]
z \hspace{10pt} & \longmapsto & \prod (a_i - b_i ) \; z
\end{array} $$
(with the convention that all terms containing a square are zero).

For all orderings of the variables such that 
$a_1 < b_1 \leq a_2 < b_2 \leq \cdots \leq a_r < b_r < x$ for all remaining variables $x$, 
$\pi $ is the identity map.
\end{theorem}

\begin{proof}
If $r=0$, let $\varphi = id $.  By induction on $r$ we reduce to the case $r=1$, since a 
composition of quasi-isomorphisms is a quasi-isomorphism and the degrees add.  To 
simplify the notation, let $a=a_1$ and $b=b_1$.  Choose an order of the variables such 
that $a,b$ are the first two variables and $a<b$.

Suppose first that $a b \in I$.  Define $\varphi $ to be the multiplication by $(a-b)$.  
Since $\Sigma _{a,b} \bigl( R(I:a) \bigr)$ is a subcomplex of $R(I)$, $\varphi (z)$ is a 
chain and we see immediately 
that it is a cycle.  The fact that $\varphi $ is a quasi-isomorphism of degree one follows 
from the Mayer-Vietoris sequence associated to the decomposition 
$\Sigma _{a,b} \bigl( R(I:a) \bigr) = A_a \bigl( R(I:a) \bigr) \cup A_b \bigl( R(I:a) \bigr)$: 
in fact ${\rm \tilde H}_* \bigl( A_a \bigl( R(I:a) \bigr) , \mathbb{Z} \bigr) \oplus 
{\rm \tilde H}_* \bigl( A_b \bigl( R(I:a) \bigr) , \mathbb{Z} \bigr) = (0)$, and 
$$ \xymatrix { 0 \ar[r] & 
{\rm \tilde H}_* \bigl( \Sigma _{a,b} \bigl( R(I:a) \bigr) , \mathbb{Z} \bigr) \ar[r] ^{\hspace{5pt} \delta } & 
{\rm \tilde H}_{*-1} \bigl( R(I:a) , \mathbb{Z} \bigr) \ar[r] & 0 } $$
is exact.  By Theorem~\ref{paba} and Remark~\ref{cretino}, the 
inclusion $\Sigma _{a,b} \bigl( R(I) \bigr) \subset R(I)$ induces an 
isomorphism $\iota $ in homology.  We are done since $\varphi $ induces in homology 
the composition $\iota \circ \delta ^{-1}$.

If $a b \notin I$, consider the ideal $I' = \bigl( I, a b' , (b')^2 \bigr) \subset \Rn [b']$.  By 
the previous case we know that the multiplication by $(a-b')$ induces an isomorphism 
of degree one between the homology of $R(I':a)$ and the homology of 
$\Sigma _{a,b'} \bigl( R(I':a) \bigr)$.  We have the following commutative diagram 
$$ \xymatrix @C=50pt {
Z \bigl( R(I':a) \bigr) \ar@{=}[d] \ar[r] ^{(a-b') \cdot \hspace{5pt}} & 
Z \bigl( \Sigma _{a,b'} \bigl( R(I':a) \bigr) \bigr) \ar@{^(->}[r] & Z \bigl( R(I') \bigr) \ar[d] ^{\bar \varphi } \\
Z \bigl( R(I:a) \bigr) \ar[rr] ^{\varphi }
& & Z \bigl( R(I) \bigr) } $$
where $\bar \varphi $ is the map of Lemma~\ref{separati}.  Since $\varphi $ is a composition 
of quasi-isomorphisms by the previous step and Lemmas~\ref{elcontre} and~\ref{separati}, 
$\varphi $ is also a quasi-isomorphism and we are done.
\end{proof}

We shall see that for the independence and dominance complexes of forests the 
map $\pi $ of Theorem~\ref{omolomap} is always the identity.

For the reader's convenience, we state explicitly the following easy results, that will be 
used frequently in the sequel.

\begin{lemma} \label{cicci}
Let $a_1, a_2 ,a_3 \in \Rn$ be distinct variables.
\begin{enumerate}
\item If $a_1$ dominates $a_2$ in $I$, then either $a_1$ dominates $a_2$ in $(I:a_3)$ 
or $R(I:a_3)$ is a cone with apex $a_2$. \label{udt0}
\item If $a_1$ dominates $a_2$ and $a_2$ 
dominates $a_3$ in $I$, then $a_1$ dominates $a_3$ in $I$. \label{udt1}
\item If $(a_1,a_2)$ is a resolution of $I$, $R(I:a_1a_2)$ is not a cone and 
$a_2$ dominates $a_3$ in $I$, then 
$(a_2 , a_1)$ is also a resolution of $I$. \label{udt2}
\end{enumerate}
\end{lemma}

\begin{proof}
(\ref{udt0}) 
Let $M$ (resp.~$M'$) be the set of minimal square-free generators of $I$ (resp.~$(I:a_3)$) 
that are divisible by $a_2$; note that $M' \subset M:a_3$.  Since $a_1$ dominates 
$a_2$ in $I$, it follows that all monomials of $M$ are divisible also by $a_1$; because 
$a_3 \neq a_1,a_2$ also all monomials of $M'$ are divisible also by $a_1$.  Moreover 
$M \neq \emptyset $.  If $M' \neq \emptyset $, then $a_1$ dominates $a_2$ in $(I:a_3)$; 
if $M' = \emptyset $, then $R(I:a_3)$ is a cone with apex $a_2$.

\noindent
(\ref{udt1}) Since $a_1$ dominates $a_2$, every square-free minimal generator of $I$ divisible by 
$a_2$ is also divisible by $a_1$.  Since $a_2$ dominates $a_3$ every square-free 
minimal generator of $I$ divisible by $a_3$ is also divisible by $a_2$ and hence by $a_1 a_2$.  

\noindent
(\ref{udt2})  If $a_1$ is the apex of a cone in $R(I)$, then everything is clear.  Otherwise we 
may assume that $a_1$ dominates $b_1$ in $I$ with $b_1 \neq a_2$: 
indeed if $a_1$ dominates $a_2$, since $a_2$ dominates $a_3$ in $I$, it follows that 
$a_1$ dominates also $a_3 \neq a_2$ in $I$ by part (\ref{udt1}).  Applying part (\ref{udt0}), 
either $a_1$ dominates $b_1$ in $(I:a_2)$ and we are done, or $R(I:a_2)$ is a cone 
with apex $b_1$ and hence also $R(I:a_2a_1)$ is a cone with apex $b_1$, contradicting 
the assumption.
\end{proof}

The symmetric group on $X$ acts naturally on $\Rn$.  This action need not 
preserve the ideal $I$, but in some cases it does.  We denote by 
$\sigma _{x y}$ the transposition of the variables $x$ and $y$.

\begin{lemma} \label{perle}
Let $A = (a_1 , \ldots , a_r)$ be a spherical resolution of $I$ and 
suppose that $R \bigl( c(A) \bigr)$ is not a cone.  Let $a \notin \{a_1 , \ldots , a_r\}$, 
$a$ dominate $b$ in $I$, and $I_i := (I:a_1 \cdots a_{i-1})$ for $i \in [r]$.  Then at least 
one of the following happens: 
\begin{enumerate}
\item the sequence $(a_1 , \ldots , a_r , a)$ is a spherical resolution; \label{aggi}
\item there is an index $i \in [r]$ such that $a$ dominates $a_i$ in $I_i$, $a_i$ dominates $a$ in 
$I_i$, the sequence 
$A' = (a_1 , \ldots , a_{i-1} , a , a_{i+1} , \ldots , a_r)$ is a spherical resolution, and 
$c(A) = \sigma _{a a_i} \bigl( c(A') \bigr)$. \label{scambio}
\end{enumerate}
\end{lemma}

\begin{proof}
For $i \in [r]$, let $a_i$ dominate $b_i$ in $I_i$.  If $b \neq a_1$, then by Lemma~\ref{cicci} 
part (\ref{udt0}) $a$ dominates $b$ in $I_1$, since $R(I_1)$ cannot be a cone because 
$R(c(A)$ is not a cone.  Thus, if $b \neq a_i$ for all $i \in [r]$, we prove that we are in 
case (\ref{aggi}) by iterating this argument.  
Otherwise, suppose that there is an index $i$ such that $b = a_i$.  Then $a_i$ 
dominates $b_i$ in $I_i$, and $a$ dominates $b = a_i$ in $I$ and hence in $I_i$ 
by Lemma~\ref{cicci} part (\ref{udt0}).  
If $a \neq b_i$, by Lemma~\ref{cicci} part (\ref{udt1}) $a$ 
dominates also $b_i$, we may replace $b$ by $b_i$ and reduce to the case $b \neq a_i$.

It remains to treat the case when $a = b_i$ and $b = a_i$: in this case $a$ and $a_i$ 
mutually dominate each other in $I_i$.  Hence every square-free minimal generator of $I_i$ 
divisible by $a$ is divisible by $a_i$ and conversely.  Thus exchanging $a$ and $a_i$ 
is an isomorphism of $\Rn$ that fixes $I_i$.  Hence the sequence 
$(a_1 , \ldots , a_{i-1} , a , a_{i+1} , \ldots , a_r)$ is a resolution of $I$.
This is case (\ref{scambio}).
\end{proof}

\begin{remark} \label{osserva}
The proof of Lemma~\ref{perle} implies that if $b \notin \{ a_1 , \ldots , a_r \}$, then 
(\ref{aggi}) certainly holds.
\end{remark}

\begin{cor} \label{orco}
Let $A = (a_1 , \ldots , a_r)$ be a spherical resolution of $I$ and let 
$a \notin \{a_1 , \ldots , a_r\}$ dominate $b$ in $I$.  Suppose that 
$R \bigl( c(A) \bigr)$ is not a cone.
\begin{enumerate}
\item If $b \notin \{a_1 , \ldots , a_r\}$, then the sequence 
$(a_1 , \ldots , a_i , a , a_{i+1} , \ldots , a_r)$ is a spherical resolution, for $0 \leq i \leq r$. 
\label{padre}
\item If $A$ is a maximal resolution, then there exists $i \in [r]$ such that $a$ dominates 
$a_i$ in $I_i$, 
$a_i$ dominates $a$ in $I_i$, $A' = (a , a_1 , \ldots , a_{i-1} , a_{i+1} , \ldots , a_r)$ is a maximal 
resolution and $c(A) = \sigma _{a a_i} \bigl( c(A') \bigr)$. \label{stredre}
\end{enumerate}
\end{cor}

\begin{proof}
(\ref{padre})  The result follows from Lemma~\ref{cicci} part (\ref{udt2}) by induction on 
$r-i$, the case $r-i=0$ being Lemma~\ref{perle} thanks to Remark~\ref{osserva}.

\noindent
(\ref{stredre})  Since $A$ is maximal, we are in case (\ref{scambio}) of Lemma~\ref{perle}.  
Thus there is an index $i \in [r]$ such that $a$ dominates $a_i$, $a_i$ dominates $a$ in $I_i$, 
and the sequence $A'' = (a_1 , \ldots , a_{i-1} , a , a_{i+1} , \ldots , a_r)$ is a spherical resolution.  From 
Lemma~\ref{perle} it follows that $c(A) = \sigma _{a a_i} \bigl( c(A'') \bigr)$.  Hence $A''$ 
is maximal.  Otherwise if $a''$ 
dominated $b''$ in $c(A'')$, then $\sigma _{a a_i} (a'')$ would dominate $\sigma _{a a_i} (b'')$ 
in $\sigma _{a a_i} \bigl( c(A'') \bigr)$.  By applying repeatedly Lemma~\ref{cicci} part (\ref{udt2}) 
we deduce that $A'$ is also a maximal resolution and it is spherical.  Clearly $c(A') = c(A'')$.
\end{proof}

We now prove the main result of this section.

\begin{theorem} \label{coco}
Let $I$ be a spherical ideal.  Then all resolutions of $I$ are spherical, all maximal 
resolutions of $I$ have the same depth, and $I$ has a unique core up to permutation 
of the variables.
\end{theorem}

\begin{proof}
Let $A = (a_1 , \ldots , a_r)$ and $A' = (a_1' , \ldots , a_s')$ be maximal resolutions of 
$I$ and suppose that $A$ is spherical.  We 
proceed by induction on $r$.  If $r=0$, then the only maximal 
resolution of $I$ is the empty resolution; thus $I$ is its own core and we are done.

Suppose that $r \geq 1$.  Note first that $a_1'$ is not the apex of a cone, since otherwise 
every maximal resolution of $I$ would contain $a_1'$ and $I$ would not admit maximal 
resolutions that are spherical.  If $a_1' \notin \{a_1, \ldots , a_r \}$, then we may apply 
Lemma~\ref{perle} to deduce that there is an index $i$ such that $a_1'$ dominates $a_i$ 
and $\bar A = (a_1 , \ldots , a_{i-1} , a_1' , a_{i+1} , \ldots , a_r)$ is a spherical resolution.  
Moreover, by Corollary~\ref{orco}, we have $c(A) = \sigma _{a_1' a_i} \bigl( c(\bar A) \bigr)$ 
and hence also $\bar A$ is maximal.  Thus we may replace $A$ by $\bar A$ and assume that 
there is an index $i$ such that $a_1' = a_i$.  Applying Corollary~\ref{orco} we may 
assume that $A = (a_1' , a_2 , \ldots , a_r)$, since changing the order of the elements 
of a maximal resolution does not affect the core of the resolution.  In this case, we 
have that $(a_2 , \ldots , a_r)$ and $(a_2' , \ldots , a_s')$ are both maximal resolutions 
of $(I:a_1')$ and the first one is spherical.  By induction we deduce that $r=s$ and that 
there exists a permutation $\sigma $ of the variables (different from $a_1'$) such that 
$$ c(A) = \bigl( (I:a_1') : a_2 \cdots a_r \bigr) = 
\sigma \Bigl( \bigl( (I:a_1') : a_2' \cdots a_r' \bigr) \Bigr) = \sigma \bigl( c(A') \bigr) $$
and $(a_2' , \ldots , a_r')$ is a spherical resolution of $(I:a_1')$.  Hence 
$A' = (a_1' , a_2' , \ldots , a_r')$ is a spherical resolution of $I$, and the proof is complete.
\end{proof}

\begin{example}
Consider the ideal $I \subset \mathbb{Z} [x_1 , \ldots , x_7]$ generated by $x_1 ^2 , \ldots , x_7^2$ 
and by the six monomials $ x_1 x_2$, $x_3 x_7$, $x_5 x_6$, $x_5 x_7$, $x_1 x_3 x_4$, and 
$x_2 x_3 x_4$.  The sequences $A_1 = (x_5,x_4)$ and $A_2 = (x_3,x_6)$ are both maximal 
spherical resolutions of $I$.  The cores are 
$$ \begin{array}{rcl}
c(A_1) & = & \bigl( x_4, \;x_5, \;x_6, \;x_7, \;x_1x_2, \;x_1x_3, \;x_2x_3 , \; 
x_1^2 , \; x_2^2 , \; x_3^2 \bigr) , \\[5pt]
c(A_2) & = & \bigl( x_3, \;x_5, \;x_6, \;x_7, \;x_1x_2, \;x_1x_4, \;x_2x_4 , \; 
x_1^2 , \; x_2^2 , \; x_4^2 \bigr) ,
\end{array} $$
and $\sigma _{x_3x_4} \bigl( c(A_1) \bigr) = c(A_2)$.
\end{example}

\section{Covers and Euler characteristics} \label{sequa}

In this section we give a combinatorial way to compute the Hilbert series of a monomial ideal. Thus we obtain the face polynomial and the Euler 
characteristic of a simplicial complex $\Delta$ in terms of covers of the Stanley-Reisner ideal of $\Delta$.  Since we will study complexes that 
are {\it a priori} known to be either contractible or homotopic to a sphere, 
Theorem~\ref{cara} may be useful to understand their homotopy type.  Note that also 
the constrictive simplicial complexes introduced in~\cite{EH} are known 
to be either contractible of homotopic to spheres.

Let $M$ be a finite set of monomials of $\Rn$.

\begin{definition}
Let $p \in \Rn$ be a monomial; a {\it cover $S$ of $p$ by $M$}, denoted by $S \rhd _M p$, 
is a subset $S \subset M$ such that $p = {\rm lcm} \{ s, s \in S \}$.  The covering polynomial 
$C_M (\ics) \in \Rn$ of $M$ is the polynomial
\begin{eqnarray*}
C_M (\ics) & := & \hspace {-10pt} \sum _{p \text{ \rm a monomial}} 
\hspace {-10pt} c_M (p) p \\[7pt]
c_M (p) & := & \sum _{S \rhd _M p} (-1)^{|S|}
\end{eqnarray*}
\end{definition}

The covering polynomial is indeed a polynomial: if $p \nmid {\rm lcm} (M)$ then there 
are no covers of $p$ and hence $c_M (p) = 0$. More precisely, if $M=\{m_1,\ldots,m_r\}$, 
then $C_M (\ics)=(1-m_1)\star \cdots \star (1-m_r)$, where $\star: \Rn \to \Rn$ is the 
$\mathbb Z-$linear distributive map defined on monomials $m,m'$ by 
$m \star m'= \lcm \{m,m'\}$. 

\begin{theorem} \label{velasco}
Let $M$ be a finite set of monomials and let $J$ be the ideal generated by $M$.  Then 
\begin{eqnarray*}
H \bigl( \Rn/J ; \ics \bigr) & = & \frac {C_M (\ics )} {\prod _{s=1} ^n (1-x_s)} .
\end{eqnarray*}
\end{theorem}

\begin{proof}
Let $r = |M|$.  The following well-known construction is called {\it Taylor's resolution}~\cite{T}.  
Let $F_0 := \Rn $ and, for $k \in [r]$, define 
$$ F_k := \bigoplus _{S \subset M, |S|=k} \Rn(-\lcm \{ S \}) . $$
Taylor's Theorem says that the complex 
$$ 0 \to F_r \to F_{r-1} \to \cdots \to F_0 \to \Rn/J \to 0 $$
with differential given by the simplicial boundary map is a (not necessarily 
minimal) free resolution of $\Rn/J$.  Recall that 
$$ H \bigl( \Rn (-a) \bigr) = \frac {a}{\prod _{s=1} ^n (1-x_s)} $$
and that the multi-graded Hilbert series is additive on exact sequences.  We find the 
following formula:
\begin{eqnarray*}
H \bigl( \Rn/J ; \ics \bigr) & = & 
\sum _{k=0} ^r (-1)^k 
\sum _{S \subset M, |S|=k} \frac {\lcm \{S\}} {\prod _{s=1} ^n (1-x_s)} = \\[10pt]
& = & \frac {\sum _{S \subset M} (-1)^{|S|} \lcm \{S\}} {\prod _{s=1} ^n (1-x_s)} = \\[10pt]
& = & \frac {C_M (\ics )} {\prod _{s=1} ^n (1-x_s)} 
\end{eqnarray*}
and the proof is complete.
\end{proof}

Let $\Delta$ be a simplicial complex and $B$ the set of minimal square-free generators of its Stanley-Reisner ideal. Let $I=(B,\icss)$.

\begin{theorem} \label{cara}
The reduced Euler characteristic of $\Delta $ is 
$$ \tilde e (\Delta ) = (-1)^{n-1} c_B (x_1 \cdots x_n) . $$
\end{theorem}

\begin{proof}
Let $X^2 := \{ \icss \}$.  By Theorem~\ref{velasco}, we deduce that 
$$ 2^n \tilde e (\Delta) = - \sum _{K \subset B \cup X^2} 
(-1)^{|K| + \deg \bigl( \lcm (K) \bigr)} . $$
Since $B \cap X^2 = \emptyset $, we have 
\begin{eqnarray*}
2^n \tilde e (\Delta) & = & 
- \sum _{S \subset B} (-1)^{|S|} \sum _{T \subset X^2} (-1)^{|T| + \deg \bigl( \lcm (S,T) \bigr)} = \\
& = & - \sum _{S \subset B} (-1)^{|S|} \sum _{\begin{array}{c}
\scriptstyle T_1 \cup T_2 \subset X^2 \\[-2pt]
\scriptstyle \lcm (T_1)|(\lcm (S))^2 \\[-2pt]
\scriptstyle (\lcm (S), \lcm (T_2))=1
\end{array}} 
(-1)^{|T_1| + |T_2| + \deg \bigl( \lcm (S,T_1,T_2) \bigr)} .
\end{eqnarray*}
Since $\deg \bigl( \lcm (S,T_1,T_2) \bigr) = |T_1|+2|T_2| + \deg \bigl( \lcm (S) \bigr)$, 
we obtain 
\begin{eqnarray*}
2^n \tilde e (\Delta) & = & 
- \sum _{S \subset B} (-1)^{|S| + \deg \bigl( \lcm (S) \bigr)} 
\sum _{\begin{array}{c}
\scriptstyle T_1 \cup T_2 \subset X^2 \\[-2pt]
\scriptstyle \lcm (T_1)|(\lcm (S))^2 \\[-2pt]
\scriptstyle (\lcm (S), \lcm (T_2))=1
\end{array}} 
(-1)^{|T_2|} .
\end{eqnarray*}
The last sum is zero unless $\lcm (S) = x_1 \cdots x_n$, and thus 
$$ 2^n \tilde e (\Delta) = 2^n (-1)^{n-1} \sum _{S \rhd _B x_1 \cdots x_n} (-1)^{|S|} $$
and we are done.
\end{proof}

In particular, 
$\tilde e (\Delta ) \equiv c_B (x_1 \cdots x_n) \equiv 
\bigl| \bigl\{ S \subset B \bigl| S \rhd _B x_1 \cdots x_n \bigr\} \bigr| \pmod 2$.

\begin{example}
Let $I = (x_1 , x_2 , x_3)^2 \subset \mathbb{Z}[x_1 , x_2 , x_3]$.  Thus 
$R(I)$ is the disjoint union of three points, $\tilde e \bigl( R(I) \bigr) = 2$, and the 
set of minimal square-free generators of $I$ is 
$B = \{ x_1 x_2 , x_1 x_3 , x_2 x_3 \}$.  There are four covers of $x_1 x_2 x_3$ by 
$B$ and we have $c_B (x_1 x_2 x_3) = 2$, as predicted by Theorem~\ref{cara}.  
We shall see in the next section that 
$I$ can be obtained as the edge ideal of a graph (the cycle with three vertices).
\end{example}

\section{The independence complex} \label{secin}

In this section we apply the techniques developed in Section~\ref{setre} 
to the independence complex of a forest.

Let $G = (V,E)$ be a graph with vertex set $V = \{ \ics \}$.  Let $\overline {G} \subset \Rn$ 
be the ideal generated by $\icss$ and by $x_i x_j$ for all $\{ x_i , x_j \} \in E$.  
The ideal $\overline {G}$ is called the {\it edge ideal of $G$} and the simplicial complex 
$R \bigl( \overline {G} \bigr)$ is called the {\it independence complex of $G$}.  
In particular, if $n = 0$, then $\overline {G} : =  (0) \subset \mathbb{Z}$.  We have 
\begin{equation} \label{proide}
\begin{array}{rcl}
\bigl( \overline {G} , v \bigr) & = & \overline {G \setminus \{v\}} \\[5pt]
\bigl( \overline {G} : v \bigr) & = & \overline {G \setminus N[v]}
\end{array} 
\end{equation}
where 
for all $S \subset V$, $G \setminus S$ is the graph obtained by removing from $G$ 
the vertices in $S$ and all the edges having a vertex in $S$ as an endpoint.

Let $a$ and $b$ be vertices of $G$; it follows from (\ref{proide}) that $a$ dominates $b$ 
if and only if $b$ is a leaf and $a$ is adjacent to $b$; the simplicial complex 
$R \bigl( \overline {G} \bigr)$ is a cone of apex $a$ if $a$ is an isolated vertex of $G$.  Thus, thanks to Theorem~\ref{collacolla}, we have proved the following.

\begin{prop} \label{semplice}
Let $F$ be a forest. Then $\overline {F}$ is simple and
\begin{itemize}
\item if $\overline F$ is conical, then $R \bigl( \overline F \bigr)$ collapses 
onto a point;
\item if $\overline F$ is spherical, then $R \bigl( \overline F \bigr)$ collapses 
onto the boundary of a cross-polytope of dimension $d \bigl( \overline F \bigr)$.\hfill $\Box$
\end{itemize}
\end{prop}
Let $G = (V,E)$ be a graph.  We let 
$$ cov_G(t) := \sum _{S {\text{ edge cover}}} t^{|S|} $$
be the generating function for the edge covers of $G$.

The next corollary gives several characterizations of when $\overline {F}$ is conical and consequently spherical.

\begin{cor} \label{facco}
Let $F$ be a forest.  The following are equivalent:
\begin{enumerate}
\item the ideal $\overline {F}$ is conical;
\item the complex $R \bigl( \overline {F} \bigr)$ is contractible; 
\item the reduced Euler characteristic $\tilde e \bigl( \overline {F} \bigr)$ is even; 
\item the number $cov_F(-1)$ is zero; 
\item the number $cov_F(1)$ of edge covers of $F$ is even; 
\item there is a sequence $\bigl( a_1 , \ldots , a_r \bigr)$ of vertices 
such that $a_i$ is adjacent to a leaf of $F_i := F \setminus N[\{a_1 , \ldots , a_{i-1}\}]$ 
and $F_{r+1}$ contains an isolated vertex; 
\item for all maximal sequences $\bigl( a_1 , \ldots , a_r \bigr)$ of vertices 
such that $a_i$ is adjacent to a leaf of $F_i := F \setminus N[\{a_1 , \ldots , a_{i-1}\}]$ 
there is an $i$ such that $F_i$ contains an isolated vertex; 
\item there is a vertex $v$ such that $R \bigl( \overline {F \setminus \{v\}} \bigr)$ and 
$R \bigl( \overline {F \setminus N[v]} \bigr)$ are either both contractible or both not 
contractible; 
\item for all vertices $v$ the complexes $R \bigl( \overline {F \setminus \{v\}} \bigr)$ and 
$R \bigl( \overline {F \setminus N[v]} \bigr)$ are either both contractible or both not 
contractible.
\end{enumerate}
\end{cor}

\begin{proof}
The equivalence of (1)-(7) follows at once from Theorem~\ref{coco}, Remark~\ref{sfeco}, 
and Theorem~\ref{cara}.  By Lemma~\ref{precontra}, we have 
$$ \tilde e (\overline {F}) = \tilde e (\overline {F} , v) - \tilde e (\overline {F} : v) . $$
Thus $\tilde e (\overline {F})$ is even if and only if 
$\tilde e (\overline {F} , v) \equiv \tilde e (\overline {F} : v) \pmod 2$.  Since $(\overline {F} , v) = 
\overline {F \setminus \{v\}}$ and $(\overline {F} : v) = \overline {F \setminus N[v]}$, 
we may conclude using the equivalence of (2) and (3).
\end{proof}

We now analyze the problem of 
computing the depth of $\overline {F}$. Hence, when $\overline {F}$ is spherical, we determine the dimension of the associated sphere. We prove that the depth of $\overline {F}$ equals the independent domination number of $F$.

The dominating sets and the independent dominating sets of a graph have been 
studied by several graph theorists (see, for instance, \cite{BC}, \cite{HHS}, \cite{LW}).  
The following lemma, which will be needed in the sequel, gives a result on the 
independent dominating sets of a forest.

\begin{lemma} \label{ids}
Let $F$ be a forest with at least one edge.  There are independent dominating sets of 
$F$ of cardinality $i(F)$ containing a vertex at distance one from a leaf.
\end{lemma}

\begin{proof}
We may assume that $F$ is a tree.  
Proceed by induction on the number of edges of $F$.  If the number of edges of $F$ is 
at most three, then the result is clear.

Suppose that $F$ is as shown below, where $T$ is a tree containing the vertex $e$. 
\begin{center}
\begin{minipage}{150pt}
$$ F = \xygraph {[] !~:{@{.}} 
!{<0pt,0pt>;<20pt,0pt>:} 
{\bullet} 
!{\save +<0pt,8pt>*\txt{$\scriptstyle a$}  \restore}
[r] {\bullet} 
!{\save +<0pt,8pt>*\txt{$\scriptstyle b$}  \restore}
[r] {\bullet} 
!{\save +<0pt,8pt>*\txt{$\scriptstyle c$}  \restore}
[r] {\bullet} 
!{\save +<0pt,8pt>*\txt{$\scriptstyle d$}  \restore}
[r] {\bullet} 
!{\save +<0pt,8pt>*\txt{$\scriptstyle e$}  \restore}
[rl] - [llll] } \; \bigcup _e \; T $$

\medskip
\centerline{The graph $F$}
\end{minipage} 
\end{center}

\bigskip

\noindent
Let $D$ be an independent  dominating set of 
$F$ of cardinality $i(F)$. If $e\in D$, then by minimality $b\in D$ and we are done. 
If $e,d\notin D$, then necessarily $a,c\in D$. Hence also $\bigl( D \setminus \{ a,c \} \bigr) \cup \{b,d\}$ 
is an independent  dominating set of $F$ of cardinality $i(F)$ and we are done. 
Suppose finally that $e\notin D$ and $d\in D$.  If $b \in D$, then we are done; otherwise $a$ 
must be in $D$ and we may replace $a$ by $b$.

The case of paths follows from what we said.  Thus we assume that $F$ is not a path.  
Let $G$ be the smallest subtree of $F$ containing all vertices whose valence in $F$ is at least three 
and let $v$ be a vertex whose valence in $G$ is at most one.  
All components of the forest $F \setminus \{v\}$, except for at most one, are paths with an 
endpoint adjacent to $v$.  By what we said above, we may assume that these paths consist of 
at most three vertices.  
For $i \in [3]$, we let $s_i$ be the number of such paths with $i$ vertices.  Note that 
$s_1+s_2+s_3 \geq 2$, since $v$ has valence at least three.  For ease of presentation we 
introduce the relevant notation in the following figure; the graph $T$ is a tree and could be empty.
\begin{center}
\begin{minipage}{150pt}
$$ \xygraph {[] !~:{@{.}} 
!{<0pt,0pt>;<20pt,0pt>:} 
{\bullet} 
!{\save +<0pt,6pt>*\txt{$\scriptstyle b_{s_2}$}  \restore}
[r] {\bullet} 
!{\save +<0pt,6pt>*\txt{$\scriptstyle a_{s_2}$}  \restore}
[dr] {\bullet} 
!{\save +<0pt,-6pt>*\txt{$\scriptstyle v$}  \restore}
[l] {\bullet} 
!{\save +<0pt,-6pt>*\txt{$\scriptstyle a_1$}  \restore}
[l] {\bullet} 
!{\save +<0pt,-6pt>*\txt{$\scriptstyle b_1$}  \restore}
[rrru] {\bullet}
[r] {\bullet}
!{\save +<0pt,6pt>*\txt{$\scriptstyle c_{s_3}$}  \restore}
[r] {\bullet} 
!{\save +<0pt,6pt>*\txt{$\scriptstyle d_{s_3}$}  \restore}
[d] {\bullet} 
!{\save +<0pt,-6pt>*\txt{$\scriptstyle d_1$}  \restore}
[l] {\bullet} 
!{\save +<0pt,-6pt>*\txt{$\scriptstyle c_1$}  \restore}
[l] {\bullet} 
!{\save +<0pt,-6pt>*\txt{$\scriptstyle x$}  \restore}
[d] {\bullet} 
!{\save +<0pt,-6pt>*\txt{$\scriptstyle e_{s_1}$}  \restore}
[l] {\cdots} [l] {\bullet}
!{\save +<0pt,-6pt>*\txt{$\scriptstyle e_1$}  \restore}
[uuuur] {T} 
- [ddd] [dr] - [uull] -[l] [d] - [rrrrr] [u] - [ll] -[dl] [dl] - [ur]
[l]
!{\save +<-10pt,13pt>*\txt{$\vdots $}  \restore}
[rrrr]
!{\save +<-10pt,13pt>*\txt{$\vdots $}  \restore}
} $$

\medskip
\centerline{The graph $F$}
\end{minipage} 
\end{center}

\bigskip

\noindent
We consider two cases.

{\bf Case 1.}  There exists an independent dominating set $D$ of $F$ with $|D| = i(F)$ 
containing $v$.  If $s_1 \geq 1$, then we are done.  If $s_3 \geq 1$, 
then we are again done, since we may suppose 
that $D$ contains $c_1 , \ldots , c_{s_3}$.  Thus we assume that $s_1 = s_3 = 0$, 
$s_2 \geq 2$, and $b_1 , \ldots , b_{s_2} \in D$.  Let us consider the tree 
$F' = \bigl( F \setminus \{b_{s_2}\} \bigr) \setminus \{a_{s_2}\}$.  
We have $i(F') = i(F) - 1$.  Indeed $D \setminus \{b_{s_2} \}$ is an 
independent dominating set of $F'$ of cardinality $i(F) - 1$; conversely, if $D'$ is an 
independent dominating set of $F'$, then $D' \cup \{b_{s_2}\}$ is an independent dominating 
set of $F$.  By induction there exists an independent dominating set $D''$ of $F'$ of cardinality 
$i(F')$ containing a vertex adjacent to a leaf $l$.  The vertex $l$ is a leaf also in $F$, 
since the unique vertex of $F'$ with a different valence in $F$ is $v$ and $v$ is not a leaf in 
$F'$.  Thus $D'' \cup \{ b_{s_2} \}$ is the required independent dominating set.

{\bf Case 2.}  Every independent dominating set of $F$ of cardinality $i(F)$ 
does not contain $v$.  Let $D$ be an independent dominating set of $F$ with $|D| = i(F)$.  
If $s_2 \geq 1$, then we may assume that $D$ contains $a_1 , \ldots , a_{s_2}$ and we 
are done.  If $s_1 \geq 1$, then we may assume that $s_3 = 0$, since otherwise $D$ 
contains $c_1 , \ldots , c_{s_3}$ by minimality.  Thus either $s_1 \geq 2$ and 
$s_2 = s_3 = 0$, or $s_1 = s_2 = 0$ and $s_3 \geq 2$.

If $s_1 \geq 2$ and $s_2 = s_3 = 0$, then we consider the tree 
$F' = F \setminus \{e_{s_1}\}$, by a similar reasoning as before we have $i(F') = i(F) - 1$ 
and we conclude by the induction hypothesis.

Suppose finally that $s_1 = s_2 = 0$ and $s_3 \geq 2$.  We may assume that 
$c_1 \notin D$; hence $x,d_1 \in D$ and by minimality $c_2 , \ldots , c_{s_3} \in D$.  This 
concludes the proof.
\end{proof}

The following result gives a strict link between dominating sets of a 
forest $F$ and resolutions of the ideal $\overline {F}$.

\begin{theorem} \label{forind}
Let $F$ be a forest; then $i(F) = d \bigl( \overline {F} \bigr)$.
If $\overline {F}$ is spherical, 
then $\gamma (F) = d \bigl( \overline {F} \bigr)$.
\end{theorem}

\begin{proof}
Let $(a_1 , \ldots , a_r)$ be any resolution; clearly $\{ a_1 , \ldots , a_r \}$ is an independent 
dominating set and we deduce that $i(F) \leq d \bigl( \overline {F} \bigr)$.  Thus 
we only need to prove that $i(F) \geq d \bigl( \overline {F} \bigr)$.

Proceed by induction on $i(F)$, the base case being clear.  We may assume that $F$ has 
at least one edge.  Let $D$ be an 
independent dominating set of minimum size.  Thanks to Lemma~\ref{ids} we may assume 
that $D$ contains a vertex $a_1$ at distance one from a leaf.  We have 
$i\bigl( F \setminus N[a_1] \bigr) = i(F) - 1$ since $D \setminus \{a_1\}$ 
is an independent dominating set of $F \setminus N[a_1]$ and if $D'$ is an 
independent dominating set of $F \setminus N[a_1]$, then $D' \cup \{a_1\}$ is 
an independent dominating set of $F$.  Moreover 
$d \bigl( \overline {F} \bigr) \leq d \bigl( \overline {F} : a_1 \bigr) + 1$.  
By induction $d \bigl( \overline {F} : a_1 \bigr) \leq i\bigl( F \setminus N[a_1] \bigr)$ 
and the first equality follows.

To prove the second statement, we proceed by induction on the number of vertices of $F$.  We may assume that 
$F$ is a tree.  If $F$ is a single edge, then the result is clear since any 
dominating set of minimum size and any maximal resolution must contain one 
of the endpoints of the edge.  Thus, without loss in generality, we only consider 
dominating sets of minimum size containing all vertices adjacent to a leaf.

Suppose that $a$ dominates $b$ in $\overline {F}$ and that the distance of $b$ 
from the closest vertex of valence different from two is at least three.  
By Theorem~\ref{coco}, $a$ can be completed to a maximal resolution of depth 
$d \bigl( \overline {F} \bigr)$ and hence we have 
$$ \begin{array}{rcl}
d \bigl( \overline {F} : a \bigr) & = & d \bigl( \overline {F} \bigr) - 1 , \\[5pt]
\gamma \bigl( F \setminus N[a] \bigr) & = & \gamma (F) - 1
\end{array} $$
since any dominating set of $F$ must contain at least one of $a$ and $b$, and 
if it contains the other vertex adjacent to $a$, then we may simply ``push it away'' 
from $a$.  By induction we have 
$d \bigl( \overline {F} : a \bigr) = \gamma \bigl( F \setminus N[a] \bigr)$ and 
we conclude in this case.

Moreover, if $a$ dominates $b_1$ and $b_2$, $b_1 \neq b_2$, then 
$d \bigl( \overline {F} \bigr) = d \bigl( \overline {F \setminus \{b_2\}} \bigr)$ and 
$\gamma ( F ) = \gamma \bigl( F \setminus \{b_2\} \bigr)$ 
and we conclude by induction.

Thus we may assume that no vertex of $F$ dominates more than one vertex 
and that the distance of a leaf from a vertex of valence at least three is at most two.  
Since $\overline {F}$ is spherical, no vertex of $F$ has two leaves at distance one and two 
respectively.  Let $v$ be a leaf of the smallest tree containing all vertices of valence 
at least three of a (non-path) component of $F$.  The forest $F \setminus \{v\}$ has at 
most one component which is not a path with an endpoint adjacent to $v$.  With our 
reductions, all path components created by removing $v$ consist of exactly one edge:
\begin{center}
\begin{minipage}{150pt}
$$ \xygraph {[] !~:{@{.}} 
!{<0pt,0pt>;<20pt,0pt>:} 
{\bullet} 
!{\save +<0pt,8pt>*\txt{$\scriptstyle b_1$}  \restore}
[d] {\bullet} 
!{\save +<0pt,8pt>*\txt{$\scriptstyle b_2$}  \restore}
[d] {\bullet} 
!{\save +<0pt,-5pt>*\txt{$\scriptstyle b_s$}  \restore}
!{\save +<10pt,13pt>*\txt{$\vdots $}  \restore}
[ruu] 
{\bullet} 
!{\save +<0pt,8pt>*\txt{$\scriptstyle a_1$}  \restore}
[d] {\bullet} 
!{\save +<0pt,8pt>*\txt{$\scriptstyle a_2$}  \restore}
[d] {\bullet} 
!{\save +<0pt,-6pt>*\txt{$\scriptstyle a_s$}  \restore}
 [ru] 
{\bullet} 
!{\save +<0pt,8pt>*\txt{$\scriptstyle v$}  \restore}
[r] {\bullet} [ur] {\bullet} [dd] {\bullet} 
[rl] : [ul] : [ur] [dl] - [lll] [u] - [r] - [dr] - [dl] - [l] } $$

\medskip
\centerline{The graph $F$}
\end{minipage} 
\end{center}

\bigskip

\noindent
We have $d \bigl( \overline {F} \bigr) = d \bigl( \overline {F \setminus \{a_s , b_s\}} \bigr) + 1$ and 
$\gamma ( F ) = \gamma \bigl( F \setminus \{a_s , b_s\} \bigr) + 1$ 
and we conclude by induction.
\end{proof}

Note that by Corollary~\ref{facco} being conical or spherical can be defined in purely 
graph theoretic terms.  As a consequence of Theorem~\ref{forind}, we have proved the following graph theoretic 
result.

\begin{cor} \label{corind}
Let $F$ be a forest such that $\overline {F}$ is spherical.  Then $\gamma (F) = i(F)$. \hfill $\Box$
\end{cor}

The problem of characterizing when the domination number equals the independent 
domination number appears in \cite{LW}; the forests for which this equality holds have 
been studied in \cite{My}.

The equality stated in Corollary~\ref{corind} may be false if $\overline {F}$ is conical; 
indeed the difference $i (F) - \gamma (F)$ can be any natural number.
\begin{center}
\begin{minipage}{150pt}
$$ \xygraph {[] !~:{@{.}} 
!{<0pt,0pt>;<20pt,0pt>:} 
!{\save +<-8pt,0pt>*\txt{$\scriptstyle l_0$}  \restore}
{\bullet} [d] 
!{\save +<0pt,3pt>*\txt{$\vdots$}  \restore}
[d] {\bullet} 
!{\save +<-8pt,0pt>*\txt{$\scriptstyle l_k$}  \restore}
[ru] {\bullet} 
[rr] {\bullet} [ur] 
!{\save +<8pt,0pt>*\txt{$\scriptstyle r_0$}  \restore}
{\bullet} [d]
!{\save +<0pt,3pt>*\txt{$\vdots$}  \restore}
[d] {\bullet} 
!{\save +<8pt,0pt>*\txt{$\scriptstyle r_k$}  \restore}
[rl] - [ul] - [ur] [dl] - [ll] - [ul] [dd] - [ru] } $$
\centerline{A tree $F$ with $i (F) - \gamma (F) = k$}
\end{minipage} 
\end{center}

\begin{example}
Consider the following tree $T$.
\begin{center}
\begin{minipage}{150pt}
$$ \xygraph {[] !~:{@{.}} 
!{<0pt,0pt>;<20pt,0pt>:} 
!{\save +<-8pt,0pt>*\txt{$\scriptstyle b_1$}  \restore}
{\bullet} [r] 
{\bullet} [d] 
!{\save +<8pt,-6pt>*\txt{$\scriptstyle a_1$}  \restore}
{\bullet} [d] 
{\bullet} [l] 
{\bullet} [u] 
{\bullet} [rr] 
{\bullet} [r] 
{\bullet} [u] 
{\bullet} [u] 
!{\save +<-8pt,0pt>*\txt{$\scriptstyle a_3$}  \restore}
{\bullet} [ddr] 
!{\save +<8pt,-6pt>*\txt{$\scriptstyle a_2$}  \restore}
{\bullet} [d] 
{\bullet} [ur] 
{\bullet} [ul]
{\bullet} 
[rl] - [dd] [ur] - [lllll] [u] - [dr] - [dl] [r] - [uu] [drr] - [uu]} $$
\centerline{The tree $T$}
\end{minipage} 
\end{center}
Let us check that the sequence of vertices $(a_1 , a_2 , a_3)$ is a spherical resolution 
of $\overline {T}$.  First of all, the vertex $a_1$ dominates $b_1$ (and each of the other 
four leaves adjacent to $a_1$).  The simplicial complex $R \bigl( \overline {T} : a_1 \bigr)$ 
is the same as the simplicial complex associated to the edge ideal of the graph 
$T_1 = T \setminus N[a_1]$:
\begin{center}
\begin{minipage}{150pt}
$$ \xygraph {[] !~:{@{.}} 
!{<0pt,0pt>;<20pt,0pt>:} 
{\bullet} [u] 
{\bullet} [u] 
!{\save +<-8pt,0pt>*\txt{$\scriptstyle a_3$}  \restore}
{\bullet} [ddr] 
!{\save +<8pt,-6pt>*\txt{$\scriptstyle a_2$}  \restore}
{\bullet} [d] 
{\bullet} [ur] 
{\bullet} [ul]
!{\save +<8pt,0pt>*\txt{$\scriptstyle b_2$}  \restore}
{\bullet} 
[rl] - [dd] [ur] - [ll] - [uu]} $$
\centerline{The tree $T_1$}
\end{minipage} 
\end{center}
We have that $a_2$ dominates $b_2$ in $\overline {T_1}$ and 
$\bigl( \overline {T_1} : a_2 \bigr)$ is the edge ideal of the graph $T_2$: 
\begin{center}
\begin{minipage}{150pt}
$$ \xygraph {[] !~:{@{.}} 
!{<0pt,0pt>;<20pt,0pt>:} 
!{\save +<-8pt,0pt>*\txt{$\scriptstyle a_3$}  \restore}
{\bullet} [d] 
!{\save +<-8pt,0pt>*\txt{$\scriptstyle b_3$}  \restore}
{\bullet} 
[rl] - [u]} $$
\centerline{The tree $T_2$}
\end{minipage} 
\end{center}
Now $a_3$ dominates $b_3$ in $\overline {T_2}$ and the ideal $\bigl( \overline {T_2} : a_3 \bigr)$ 
is the edge ideal of the empty tree.  Thus 
$R \bigl( \overline {T} \bigr) \simeq \Sigma ^3 \bigl( S^{-1} \bigr) \simeq S^2$.  Note that the 
independent domination number and the domination number of $T$ are both equal to 3, as 
predicted by Corollary~\ref{corind}.
\begin{center}
\begin{minipage}{100pt}
$$ \xygraph {[] !~:{@{.}} 
!{<0pt,0pt>;<15pt,0pt>:} 
{\bullet} [r] 
{\bullet} [d] 
*\cir<4pt>{}
{\bullet} [d] 
{\bullet} [l] 
{\bullet} [u] 
{\bullet} [rr] 
{\bullet} [r] 
{\bullet} [u] 
{\bullet} [u] 
*\cir<4pt>{}
{\bullet} [ddr] 
*\cir<4pt>{}
{\bullet} [d] 
{\bullet} [ur] 
{\bullet} [ul]
{\bullet} 
[rl] - [dd] [ur] - [lllll] [u] - [dr] - [dl] [r] - [uu] [drr] - [uu]} $$
\end{minipage} 
\begin{minipage}{100pt}
$$ \xygraph {[] !~:{@{.}} 
!{<0pt,0pt>;<15pt,0pt>:} 
{\bullet} [r] 
{\bullet} [d] 
*\cir<4pt>{}
{\bullet} [d] 
{\bullet} [l] 
{\bullet} [u] 
{\bullet} [rr] 
{\bullet} [r] 
{\bullet} [u] 
*\cir<4pt>{}
{\bullet} [u] 
{\bullet} [ddr] 
*\cir<4pt>{}
{\bullet} [d] 
{\bullet} [ur] 
{\bullet} [ul]
{\bullet} 
[rl] - [dd] [ur] - [lllll] [u] - [dr] - [dl] [r] - [uu] [drr] - [uu]} $$
\end{minipage} 
\centerline{The dominating sets of minimum cardinality}
\end{center}
\end{example}

We conclude this section with a corollary of Theorem~\ref{omolomap} 
giving explicitly a generator of the reduced homology of the independence 
complex $R \bigl( \overline F \bigr)$ of a forest $F$.

\begin{cor} \label{generind}
Let $F$ be a forest such that $\overline F$ is spherical, let $A= (a_1 , \ldots , a_r)$ be 
a maximal resolution of $\overline F$ and suppose that $a_i$ dominates $b_i$ in 
$\bigl( \overline F : a_1 \cdots a_{i-1} \bigr)$.  Then there is an order of the variables 
such that $a_1 < b_1 < a_2 < b_2 < \cdots < a_r < b_r$ are the first $2r$ variables; with 
such an order, the reduced homology of $R\bigl( \overline F \bigr)$ is generated by the 
cycle 
$$ z := \prod _{i=1} ^r (a_i - b_i) . $$
\end{cor}

\begin{proof}
The core $c(A)$ is the ideal generated by $X$.  Hence a homology generator of 
$R \bigl( c(A) \bigr)$ is the class associated to the cycle 1 (in degree $-1$).  Since 
$a_i , b_i \in \bigl( \overline F : a_1 \cdots a_i \bigr)$ for all $i$, it follows that all 
the variables $a_i$ and $b_i$ are distinct.  Thus there exists an order of the variables 
such that $a_1 < b_1 < a_2 < b_2 < \cdots < a_r < b_r$ are the first $2r$ variables, 
and the result follows by Theorem~\ref{omolomap}.
\end{proof}

\section{Further results on the independence complex} \label{sesei}

In this section, we give some results on the independence complex of a general graph $G$.  
Note that the independence complex of a graph may have any homotopy type 
(see~\cite{EH}, Section~9).  We give some restrictions to the overall 
complexity of the independence complex of $G$.

Given a topological space $T$, let $h(T)$ be the sum (if finite) of the ranks of all its reduced 
homology groups and let $hd(T)$ be $-\infty $ if $\tilde H_k (T , \mathbb{Z}) = (0)$ for all $k$ 
and let $hd (T)$ be the maximum (if finite) of the integers $k$ such that 
$\tilde H_k (T , \mathbb{Z}) \neq (0)$ otherwise.  Both $h(T)$ and $hd(T)$ are rough measures 
of how complicated $T$ is.  
In particular, suppose that $T$ is a simply connected topological space; it is known that 
$T$ is contractible or homotopy equivalent to a sphere if and only if $h(T) \leq 1$.  Let 
$\Delta $ be a simplicial complex different from $S^{-1}$.  We denote by $h(\Delta )$ and 
$hd(\Delta )$ the corresponding functions on the realization of $\Delta $.  Observe that 
$h(S^{-1}) = 1$ and $hd(S^{-1}) = -1$.

Let $G = (V,E)$ be a finite graph with $\kappa $ connected components.  
Let $h_1(G) := \kappa + |E| - |V|$, the rank of the first homology group of the 
topological space underlying $G$.  The following Theorem gives sharp upper 
bounds for $h \bigl( R \bigl( \overline {G} \bigr) \bigr)$ and 
$hd \bigl( R \bigl( \overline {G} \bigr) \bigr)$ in terms of $h_1(G)$ and $|V|$ 
respectively.

\begin{theorem}
Let $G$ be a finite graph; then $h \bigl( R \bigl( \overline {G} \bigr) \bigr) \leq 2^{h_1(G)}$ 
and $hd \bigl( R \bigl( \overline {G} \bigr) \bigr) \leq |V|/2-1$.
\end{theorem}

\begin{proof}
We prove first that $h \bigl( R \bigl( \overline {G} \bigr) \bigr) \leq 2^{h_1(G)}$.  Proceed 
by induction on $h_1(G)$.  If $h_1(G) = 0$, then $G$ is a forest and the result follows 
from Proposition~\ref{semplice}.  Suppose that $h_1(G) \geq 1$ and let $v$ be a vertex of 
$G$ such that $h_1 \bigl( G \setminus \{v\} \bigr) < h_1(G)$ (it suffices to choose a 
vertex $v$ contained in a cycle of $G$).  By Lemma~\ref{precontra} and the 
Mayer-Vietoris sequence we have 
$$ h \Bigl( R \bigl( \overline {G} \bigr) \Bigr) \leq h \Bigl( R \bigl( \overline {G \setminus \{v\}} \bigr) \Bigr) + 
h \Bigl( R \bigl( \overline {G \setminus N[v]} \bigr) \Bigr) . $$
Thus from the inductive hypothesis and the choice of $v$ we conclude that 
$$ h \Bigl( R \bigl( \overline {G} \bigr) \Bigr) \leq 
2^{h_1 \left( \vphantom{M^M} G \setminus \{v\} \right)} + 
2^{h_1 \left( \vphantom{M^M} G \setminus N[v] \right)} \leq 
2^{h_1 ( G )} . $$
We now prove that $hd \bigl( R \bigl( \overline {G} \bigr) \bigr) \leq |V|/2-1$.  Proceed by induction on 
$|V|$.  If $|V| = 0$, then the assertion is clear.  Suppose that $|V| \geq 1$ and 
let $v$ be a vertex of $G$.  If $v$ is an isolated vertex, then $R \bigl( \overline {G} \bigr)$ is contractible 
and the result follows.  If $v$ has valence at least one, by 
Lemma~\ref{precontra} and the inductive hypothesis we have 
\begin{eqnarray*}
hd \Bigl( R \bigl( \overline {G} \bigr) \Bigr) & \leq & 
\max \Bigl\{ hd \Bigl( R \bigl( \overline {G \setminus \{v\}} \bigr) \Bigr) , 
hd \Bigl( R \bigl( \overline {G \setminus N[v]} \bigr) \Bigr) + 1 \Bigr\} \leq \\[5pt]
& \leq & \max \Bigl\{ \frac {|V| -1} {2} - 1 , \frac {|V| - 2} {2} \Bigr\} = 
\frac {|V|} {2} - 1 
\end{eqnarray*}
and this concludes the proof.
\end{proof}

The bounds are sharp: for all $n \geq 1$, the disjoint union $T^{(n)}$ of $n$ triangles realizes the 
first upper bound and the disjoint union $P^{(n)}$ of $n$ edges realizes the second one.  Indeed 
$$ h \bigl( T^{(n)} \bigr) \geq \Bigl| \tilde e \Bigl( R \Bigl( \overline {T^{(n)}} \Bigr) \Bigr) \Bigr| = 
\tilde e \Bigl( R \Bigl( \overline {T^{(1)}} \Bigr) \Bigr)^n = 2^n $$
thanks to Theorem~\ref{cara}, and 
$$ hd \Bigl( R \Bigl( \overline {P^{(n)}} \Bigr) \Bigr) = 
d \Bigl( \overline {P^{(n)}} \Bigr) - 1 = n - 1 $$
thanks to Theorem~\ref{forind}, since $\overline {P^{(n)}}$ is spherical 
by Corollary~\ref{facco}.

In the remainder of this section we study the independence complex of $G$, when 
$h_1 (G) = 1$.  In particular we determine the homotopy type of 
$R \bigl( \overline {G} \bigr)$ and we establish when $\overline {G}$ is simple.

\begin{lemma} \label{ciclo}
Let $C_k$ be the cycle with $k$ vertices.  We have 
$$ \begin{array}{lcll}
\overline {C_{3n-1}} & \simeq & S^{n-1} \\[5pt]
\overline {C_{3n}} & \simeq & S^{n - 1} \wedge S^{n - 1} \\[5pt]
\overline {C_{3n+1}} & \simeq & S^{n - 1} . 
\end{array} $$
\end{lemma}

\begin{proof}
Let $P_l$ denote the path with $l$ vertices and let $v$ be a vertex of $C_k$.  
By Lemma~\ref{precontra} 
$$ \overline {C_k} = A_v \bigl( \overline {P_{k-3}} \bigr) 
\bigcup _{\overline {P_{k-3}}} \overline {P_{k-1}} . $$
If $k = 3n-1$, then $\overline {P_{k-1}}$ is contractible and hence 
$\overline {C_k} \simeq \Sigma \bigl( \overline {P_{k-3}} \bigr) \simeq 
\Sigma \bigl( S^{n-2} \bigr) \simeq S^{n-1}$.

\noindent
If $k = 3n$, then $\overline {P_{k-3}} \simeq S^{n-2}$ and 
$\overline {P_{k-1}} \simeq S^{n-1}$.  Thus $\overline {P_{k-3}}$ is homotopic to a point 
in $\overline {P_{k-1}}$ and $\overline {C_k} \simeq S^{n-1} \wedge S^{n-1}$.

\noindent
If $k = 3n+1$, then $\overline {P_{k-3}}$ is contractible and hence 
$\overline {C_k} \simeq \overline {P_{k-1}} \simeq S^{n-1}$.
\end{proof}

If $G$ is a graph with $h_1(G) = 1$, then $G$ contains a unique cycle; let $n$ be the 
length of the cycle and denote the cycle by $C_n$.  For each vertex $v$ on $C_n$ 
there are $t_v \geq 0$ trees $T_1 ^v , \ldots , T_{t_v} ^v$ attached to $v$ in such a 
way that a leaf $l_i^v$ of $T_i ^v$ is identified with $v$.  We call the $T_i ^v$'s the 
{\it tree tentacles of $G$}.  Note that $G$ may contain connected components which 
are trees.

\begin{prop} \label{senciclo}
Let $G$ be a graph with $h_1(G) = 1$.  Then either $\overline {G}$ is simple or 
$c(I) = \{ \overline {C_n} \}$.
\end{prop}

\begin{proof}
Let $A$ be a maximal resolution of $\overline {G}$ and note that $c(A)$ is the edge ideal of a 
subgraph $\Gamma $ of $G$.  Since $A$ is maximal, $\Gamma $ has neither isolated 
vertices, nor leaves.  It follows that $\Gamma $ is either empty or $\Gamma = C_n$, the 
unique cycle of $G$.
\end{proof}

\begin{lemma} \label{tenta}
Let $G$ be a graph with $h_1(G)=1$.  
If a tree tentacle $T$ of $G$ is such that $\overline {T}$ is spherical, then $\overline {G}$ is simple.
\end{lemma}

\begin{proof}
Suppose that $l$ is the leaf of $T$ identified with the vertex $v$ of $G$.  Proceed by 
induction on the depth of $\overline {T}$.  If $d \bigl( \overline {T} \bigr) = 1$, then $T$ 
is a path with two or three vertices and removing the vertex $a$ of $T$ adjacent to the 
leaf different from $l$ shows that 
$\overline {G} \simeq \Sigma \bigl( \overline {G \setminus N[a]} \bigr)$.  
Since $G \setminus N[a]$ is a forest, it follows that $\overline {G}$ is simple.

Suppose that $d \bigl( \overline {T} \bigr) \geq 2$ and note that $T$ has a leaf $b$ 
different from $l$.  Let $a$ be the vertex of $T$ adjacent to $b$; then $a$ dominates 
$b$ in $G$, we have $\overline {G} \simeq \Sigma \bigl( \overline {G \setminus N[a]} \bigr)$.  
By induction $\overline {G \setminus N[a]}$ is simple and we conclude that 
$\overline {G}$ is simple.
\end{proof}

\begin{theorem} \label{casi}
Let $G$ be a graph with $h_1(G)=1$.  We have the following possibilities: 
\begin{enumerate}
\item $R \bigl( \overline {G} \bigr)$ is contractible; 
\item $R \bigl( \overline {G} \bigr)$ is homotopic to a sphere; 
\item $R \bigl( \overline {G} \bigr)$ is homotopic to a wedge of two spheres of the same dimension. 
\label{wedge}
\end{enumerate}
Moreover (\ref{wedge}) happens only if for every tree tentacle $T$ of $G$ the ideal $\overline {T}$ 
is conical and $n \equiv 0 \pmod 3$, where $n$ is the length of the unique cycle of $G$.
\end{theorem}

\begin{proof}
If $\overline {G}$ is simple, then the result is clear.  Otherwise, by Theorem~\ref{collacolla} and 
Proposition~\ref{senciclo}, $R \bigl( \overline {G} \bigr)$ collapses to either a point or 
to an iterated suspension of $R \bigl( \overline {C_n} \bigr)$.  The first part of the Theorem 
follows by Lemma~\ref{ciclo}.

By Lemma~\ref{tenta} and Proposition~\ref{senciclo}, case (\ref{wedge}) cannot happen 
unless for every tree tentacle $T$ the ideal $\overline {T}$ is conical and $c(I) = \{ \overline {C_n} \}$.  
By Lemma~\ref{ciclo}, 
$R \bigl( \overline {C_n} \bigr)$ is a wedge of two spheres only if $n \equiv 0 \pmod 3$.
\end{proof}

The following result, whose proof is straightforward, states that the Euler characteristic of 
$R \bigl( \overline {G} \bigr)$ determines in which of the cases of Theorem~\ref{casi} we are.  
We note that the reduced Euler characteristic can be computed using covering numbers 
by Theorem~\ref{cara}.  Thus we may determine (up to an iterated suspension) the 
homotopy type of $R \bigl( \overline {G} \bigr)$ without constructing any resolution of $\overline {G}$.

\begin{cor}
Let $G$ be a graph with $h_1(G) = 1$ and let $\Delta := R \bigl( \overline {G} \bigr)$.
\begin{enumerate}
\item If $cov (G) = 0$, then $\Delta $ is contractible.
\item If $\bigl| cov (G) \bigr| = 1$, then $\Delta $ is homotopic to a sphere.
\item If $\bigl| cov (G) \bigr| = 2$, then $\Delta $ is homotopic to a wedge 
of two spheres of the same dimension. \hfill $\Box$
\end{enumerate} 
\end{cor}

\section{The dominance complex} \label{seset}

In this section we apply the techniques that we developed in Section~\ref{setre} 
to the dominance complex of a forest.

Let $G$ be a graph with vertices $\ics$.  Let $G^\star$ be the ideal generated by 
$\bigl\{ \prod _{x \in N[x_i]} x \bigr\} _{i=1} ^n$ and $\icss $.  The ideal 
$G^\star$ is called the {\it star ideal of $G$} and the simplicial complex 
$R \bigl( G^\star \bigr)$ is called the {\it dominance complex of $G$}.  
The dominance complex of $G$ is never a cone, since every variable divides some minimal generator of $G^\star$.

Let $a \in X$; we have 
$$ \bigl( G^\star : a \bigr) = 
\Bigl( \bigl( G \setminus \{ {\text{edges containing }} a \} \bigr) ^\star , 
\prod _{y \in N[a] \setminus \{a\}} y \Bigr) . $$
If $a$ dominates $b$, then $b$ is adjacent to $a$ and all vertices adjacent to $b$ 
are also adjacent to $a$, i.e.~$N[b] \subset N[a]$.  Hence if $a$ dominates $b$ we 
have $\bigl( G^\star : a \bigr) = \bigl( G \setminus \{a\} \bigr) ^\star $.  In particular the 
dominance complex of a graph is never a cone.  Thus, when $F$ is a forest, $F^\star$ 
is always spherical and simple, since the vertex adjacent to a leaf dominates the leaf; 
this is the unique way a vertex may dominate another one in a forest.  The following 
theorem relates the dominance complex of a forest $F$ to the matching number 
$\beta _1 (F)$ and the vertex covering number $\alpha _0 (F)$, which are known to be 
equal (see Theorem~\ref{konig}).

\begin{theorem} \label{mucca}
Let $F$ be a forest; then 
\begin{enumerate}
\item $F^\star $ is simple;  \label{stellastalla}
\item the dominance complex of $F$ collapses onto the boundary of a cross-polytope 
of dimension $d \bigl( F^\star \bigr)$; \label{croce}
\item $\beta _1 (F) = \alpha _0 (F) = d \bigl( F^\star \bigr)$. \label{alfabeta}
\end{enumerate}
\end{theorem}

\begin{proof}
(\ref{stellastalla}) Follows from the remarks preceding the statement of the theorem.

\noindent
(\ref{croce}) Follows from Theorem~\ref{collacolla} since if $a$ dominates $b$ in $I$, 
then $b$ is a leaf and $ab = N[b] \in I$.

\noindent
(\ref{alfabeta}) 
Let $\bigl( a_1 , \ldots , a_r \bigr)$ be a maximal resolution of $F$ and suppose that 
$a_i$ dominates $b_i$, for $i \in [r]$.  Clearly $\bigl\{ \{a_1,b_1\} , 
\ldots , \{a_r,b_r\} \bigr\}$ is a matching of $F$ and $\bigl\{ a_1 , \ldots , a_r \bigr\}$ 
is a vertex cover of $F$.  Thus it suffices to show that $|M| \leq r$ for all 
matchings $M$ of $F$ and that $|C| \geq r$ for all vertex covers $C$ of 
$F$.  Proceed by induction on $r$.  The result is clear when $r=0$, since in this 
case $F$ has no edges, and $\emptyset $ is both the unique matching and the minimum 
vertex cover of $F$.  Suppose that $r \geq 1$.  Let $F' := F \setminus \{a_1\}$, 
$M' \subset M$ be the matching of $F'$ induced by $M$, and $C' = C \setminus \{a_1,b_1\}$.  
Note that $|M'| \geq |M|-1$ and $|C'| \leq |C| - 1$.  
Since $\bigl( a_2 , \ldots , a_r \bigr)$ is a maximal resolution of $F'$, the assertion 
follows by the inductive hypothesis.
\end{proof}

A consequence of Theorem~\ref{mucca} is that the removal of a single vertex of $F$ 
decreases the depth of $F^\star$ by at most one: given any matching $M$ of maximum 
size of $F$, removing a vertex forces the removal of at most one edge from $M$.

\begin{example}
Consider the following tree $T$.
\begin{center}
\begin{minipage}{150pt}
$$ \xygraph {[] !~:{@{.}} 
!{<0pt,0pt>;<20pt,0pt>:} 
!{\save +<-8pt,0pt>*\txt{$\scriptstyle b_1$}  \restore}
{\bullet} [r] 
{\bullet} [d] 
!{\save +<8pt,-6pt>*\txt{$\scriptstyle a_1$}  \restore}
{\bullet} [d] 
{\bullet} [l] 
{\bullet} [u] 
{\bullet} [rr] 
{\bullet} [r] 
!{\save +<0pt,-6pt>*\txt{$\scriptstyle a_3$}  \restore}
{\bullet} [u] 
!{\save +<-8pt,0pt>*\txt{$\scriptstyle a_2$}  \restore}
{\bullet} [u] 
{\bullet} [ddr] 
!{\save +<8pt,-6pt>*\txt{$\scriptstyle a_4$}  \restore}
{\bullet} [d] 
{\bullet} [ur] 
{\bullet} [ul]
{\bullet} 
[rl] - [dd] [ur] - [lllll] [u] - [dr] - [dl] [r] - [uu] [drr] - [uu]} $$
\centerline{The tree $T$}
\end{minipage} 
\end{center}
Let us check that the sequence of vertices $(a_1 , a_2 , a_3 , a_4)$ is a resolution 
of $T ^\star$.  First of all, the vertex $a_1$ dominates $b_1$ (and each of the other 
four leaves adjacent to $a_1$).  The ideal $\bigl( T^\star : a_1 \bigr)$ is the star 
ideal of the graph 
$T_1 = T \setminus \{ {\text{edges containing }} a_1 \}$:
\begin{center}
\begin{minipage}{150pt}
$$ \xygraph {[] !~:{@{.}} 
!{<0pt,0pt>;<20pt,0pt>:} 
{\bullet} [r] 
{\bullet} [d] 
!{\save +<1pt,-6pt>*\txt{$\scriptstyle a_1$}  \restore}
{\bullet} 
[d] 
{\bullet} [l] 
{\bullet} [u] 
{\bullet} [rr] 
{\bullet} [r] 
!{\save +<0pt,-6pt>*\txt{$\scriptstyle a_3$}  \restore}
{\bullet} [u] 
!{\save +<-8pt,0pt>*\txt{$\scriptstyle a_2$}  \restore}
{\bullet} [u] 
!{\save +<-8pt,0pt>*\txt{$\scriptstyle b_2$}  \restore}
{\bullet} [ddr] 
!{\save +<8pt,-6pt>*\txt{$\scriptstyle a_4$}  \restore}
{\bullet} [d] 
{\bullet} [ur] 
{\bullet} [ul]
{\bullet} 
[rl] - [dd] [ur] - [lll] [r] - [uu]} $$
\centerline{The forest $T_1$}
\end{minipage} 
\end{center}
Note that the simplicial complexes associated to $T_1^\star$ and to the star ideal 
of the graph $T_1 \setminus \bigl\{ \text{isolated vertices of } T_1 \bigr\}$ are the 
same; thus in what follows we will always remove isolated vertices.  We have that 
$a_2$ dominates $b_2$ in $T_1 ^\star $ and $\bigl( T_1 ^\star : a_2 \bigr)$ is the 
star ideal of the graph $T_2$: 
\begin{center}
\begin{minipage}{150pt}
$$ \xygraph {[] !~:{@{.}} 
!{<0pt,0pt>;<20pt,0pt>:} 
!{\save +<0pt,-6pt>*\txt{$\scriptstyle b_3$}  \restore}
{\bullet} [r] 
!{\save +<0pt,-6pt>*\txt{$\scriptstyle a_3$}  \restore}
{\bullet} 
[r] 
!{\save +<8pt,-6pt>*\txt{$\scriptstyle a_4$}  \restore}
{\bullet} [d] 
{\bullet} [ur] 
{\bullet} [ul]
{\bullet} 
[rl] - [dd] [ur] - [lll] [r] } $$
\centerline{The tree $T_2$}
\end{minipage} 
\end{center}
Now $a_3$ dominates $b_3$ in $T_2 ^\star$ and the ideal $\bigl( T_2 ^\star : a_3 \bigr)$ is the star 
ideal of the tree $T_3$: 
\begin{center}
\begin{minipage}{150pt}
$$ \xygraph {[] !~:{@{.}} 
!{<0pt,0pt>;<20pt,0pt>:} 
!{\save +<-8pt,0pt>*\txt{$\scriptstyle b_4$}  \restore}
{\bullet} [d] 
!{\save +<-8pt,0pt>*\txt{$\scriptstyle a_4$}  \restore}
{\bullet} [d] 
{\bullet} [ur] 
{\bullet} 
[rl] - [l] [u] - [dd]} $$
\centerline{The tree $T_3$}
\end{minipage} 
\end{center}
Finally $a_4$ dominates $b_4$ in $T_3 ^\star$ and the ideal $\bigl( T_3 ^\star : a_4 \bigr)$ 
is the star ideal of the empty graph.  Thus 
$R \bigl( T^\star \bigr) \simeq \Sigma ^4 \bigl( S^{-1} \bigr) \simeq S^3$.  Note that the 
matching number and the vertex covering number of $T$ are both equal to 4, as predicted 
by Theorem~\ref{mucca}.
\begin{center}
\begin{minipage}{100pt}
$$ \xygraph {[] !~:{@{=}} 
!{<0pt,0pt>;<15pt,0pt>:} 
!{<0pt,0pt>;<15pt,0pt>:} 
{\bullet} [r] 
{\bullet} [d] 
{\bullet} [d] 
{\bullet} [l] 
{\bullet} [u] 
{\bullet} [rr] 
{\bullet} [r] 
{\bullet} [u] 
{\bullet} [u] 
{\bullet} [ddr] 
{\bullet} [d] 
{\bullet} [ur] 
{\bullet} [ul]
{\bullet} 
[rl] : [d] -[d] [ur] - [ll] : [l] - [ll] [u] : [dr] - [dl] [r] - [uu] [drr] - [u] : [u]} $$
\end{minipage} 
\begin{minipage}{100pt}
$$ \xygraph {[] !~:{@{.}} 
!{<0pt,0pt>;<15pt,0pt>:} 
{\bullet} [r] 
{\bullet} [d] 
*\cir<4pt>{}
{\bullet} [d] 
{\bullet} [l] 
{\bullet} [u] 
{\bullet} [rr] 
{\bullet} [r] 
*\cir<4pt>{}
{\bullet} [u] 
*\cir<4pt>{}
{\bullet} [u] 
{\bullet} [ddr] 
*\cir<4pt>{}
{\bullet} [d] 
{\bullet} [ur] 
{\bullet} [ul]
{\bullet} 
[rl] - [dd] [ur] - [lllll] [u] - [dr] - [dl] [r] - [uu] [drr] - [uu]} $$
\end{minipage} 
\centerline{A matching and a vertex cover of maximum cardinality}
\end{center}
\end{example}

We prove a corollary of Theorem~\ref{omolomap} entirely analogous to 
Corollary~\ref{generind}.  Note that despite the fact that the dominance 
complex is substantially different from the independence complex, 
the statements and the proofs of Corollaries~\ref{generind} and~\ref{generdom} 
are identical.

\begin{cor} \label{generdom}
Let $F$ be a forest and $A= (a_1 , \ldots , a_r)$ be a maximal resolution of $F ^\star$, 
and suppose that $a_i$ dominates $b_i$ in 
$\bigl( F ^\star : a_1 \cdots a_{i-1} \bigr)$.  Then there is an order of the variables 
such that $a_1 < b_1 < a_2 < b_2 < \cdots < a_r < b_r$ are the first $2r$ variables; with 
such an order, the reduced homology of $R\bigl( F ^\star \bigr)$ is generated by the 
cycle 
$$ z := \prod _{i=1} ^r (a_i - b_i) . $$
\end{cor}

\begin{proof}
The core $c(A)$ is the ideal generated by $X$.  Hence a homology generator of 
$R \bigl( c(A) \bigr)$ is the class associated to the cycle 1 (in degree $-1$).  Since 
$a_i , b_i \in \bigl( F ^\star : a_1 \cdots a_i \bigr)$ for all $i$, it follows that all 
the variables $a_i$ and $b_i$ are distinct.  Thus there exists an order of the variables 
such that $a_1 < b_1 < a_2 < b_2 < \cdots < a_r < b_r$ are the first $2r$ variables, 
and the result follows by Theorem~\ref{omolomap}.
\end{proof}

Clearly, if $F$ has no isolated vertices, then $\gamma (F) \leq \alpha _0 (F)$.  Thus 
$d \bigl( \overline {F} \bigr) \leq d \bigl( F^\star \bigr)$.  In particular, if $d \bigl( \overline {F} \bigr)$ 
is spherical, then the sphere associated to the independence complex has dimension at 
most the dimension of the sphere associated to the dominance complex.

We conclude with the following consequence of Theorem~\ref{mucca}.

\begin{cor}
Let $F = (V,E)$ be a forest; we have 
$$ \sum _{S {\text{ dominating set}}} \hspace{-10pt} (-1)^{|S|} = (-1) ^{\beta _1 (F) + |V|} . $$
\end{cor}

\begin{proof}
Follows at once from Theorem~\ref{cara} and Theorem~\ref{mucca}.
\end{proof}

The previous result could be also proved without Theorem~\ref{cara}. Indeed, it is enough 
to consider the map sending a subset of $V$ to its complementary since the faces of the 
dominance complex are the complements of the dominating sets.  
Note also that, if $F$ has no trivial components, then $\beta _1 (F) + |V| \equiv \alpha_1(F) \pmod 2$ 
by Theorem~\ref{gallai}.

\section{Acknowledgments}
We are indebted to Mauricio Velasco who simplified the proof of Theorem~\ref{velasco}.


\begin{thebibliography}{99}
\bibitem[AL]{AL} Allan, R.B., Laskar, R., {\it On domination and independent domination numbers 
of a graph}, Discrete Mathematics, 23 (1978), 73-76.
\bibitem[ALH]{ALH} Allan, R.B., Laskar, R., Hedetniemi S., {\it A note on total domination}, 
Discrete Mathematics, 49 (1984), 7-13.
\bibitem[B]{B} Bollob\'as, B., {\it Modern Graph Theory}, Graduate Texts in 
Mathematics, Vol.~184, Springer, 1998.
\bibitem[BC]{BC} Bollob\'as, B., Cockayne, E.J., {\it Graph-theoretic parameters concerning domination, independence, and irredundance}, J.~Graph Theory 3 (1979), no.~3, 241-249. 
\bibitem[C]{C} Cohen, M.M., {\it A course in simple-homotopy theory}, 
Graduate Texts in Mathematics, Vol.~10, Springer-Verlag, New York-Berlin, 1973.
\bibitem[D]{D} Diestel, R., {\it Graph Theory}, Graduate Texts in 
Mathematics, Vol.~173, Springer, 1997.
\bibitem[EH]{EH} Ehrenborg, R., Hetyei, G., {\it The topology of the independence complex}, 
European J.~Combin.~27 (2006), no.~6, 906-923.
\bibitem[ET]{ET} Erd\"{o}s, P., Tuza, Z., {\it Vertex coverings of the edge set in a connected graph},  Graph theory, combinatorics, and algorithms, Vol.~1, 2 (Kalamazoo, MI, 1992),  1179-1187, Wiley-Intersci.~Publ., Wiley, New York, 1995.
\bibitem[HHS]{HHS} Haynes, T.W., Hedetniemi, S.T., Slater, P.J., {\it Fundamentals of domination in graphs}, Monographs and Textbooks in Pure and Applied Mathematics, 208.  Marcel Dekker, Inc., New York, 1998.
\bibitem[HY]{HY} Henning, M.A., Yeo, A., {\it Total domination and matching numbers in claw-free graphs}, Electron.~J.~Combin.~13 (2006), no.~1, Research Paper 59, 28 pp.
\bibitem[LW]{LW} Laskar, R., Walikar, H.B., {\it On domination related concepts in graph theory}, 
Combinatorics and Graph Theory, Lecture Note in Mathematics, 885 (1981), 308-320, Springer.
%
\bibitem[MS]{MS} Miller, E., Sturmfels, B., {\it Combinatorial commutative algebra}, 
Graduate Texts in Mathematics, Vol.~227, Springer-Verlag, New York, 2005.
%
\bibitem[Mu]{Mu} Munkres, J.R., {\it Elements of Algebraic Topology},  Perseus Books Publishing, 1984.
\bibitem[My]{My} Mynhardt, C.M., {\it Vertices contained in every minimum dominating set of a tree},  J.~Graph Theory  31  (1999),  no.~3, 163-177.
\bibitem[S]{S} Stanley, R., {\it Combinatorics and Commutative Algebra}, Second edition, 
Progress in Mathematics, Vol.~41, Birkh\"auser, 1996.
\bibitem[T]{T} Taylor, D., {\it Ideals generated by monomials in an $R-$sequence}, PhD Thesis, 
Univ.~of Chicago, 1960.
\end{thebibliography}
\end{document}